\documentclass[11pt]{article}
\usepackage{amsthm,amssymb,amsmath}
\usepackage{epsfig}
\usepackage{color}

\newtheorem{prelem}{{\bf Theorem}}

\newtheorem{theorem}{Theorem}
\newtheorem{corollary}[theorem]{Corollary}

\newtheorem{lemma}[theorem]{Lemma}

\newtheorem{conj}[theorem]{Conjecture}
\theoremstyle{definition}

\theoremstyle{remark}

\begin{document}

\title{Roman domination in graphs with minimum degree at least
two and some forbidden cycles }
\date{}
\author{ S.M. Sheikholeslami$^{1}$\thanks{%
Corresponding author}, M. Chellali$^2$, R. Khoeilar$^1$, H. Karami$^1$ and
Z. Shao$^3$, \vspace{2mm} \\
$^1$Department of Mathematics \\
Azarbaijan Shahid Madani University\\
Tabriz, I.R. Iran\\
\texttt{s.m.sheikholeslami@azaruniv.ac.ir}\\
\texttt{khoeilar@azaruniv.ac.ir}\\
\texttt{h.karami@azaruniv.ac.ir}\vspace{3mm}\\
$^2$LAMDA-RO Laboratory, Department of Mathematics\\
University of Blida\\
B.P. 270, Blida, Algeria\\
\texttt{m\_chellali@yahoo.com}\vspace{2mm}\\
$^3$Institute of Computing Science and Technology\\
Guangzhou University\\
Guangzhou 510006, China\\
\texttt{zshao@gzhu.edu.cn}\vspace{2mm}\\
}
\maketitle

\begin{abstract}
Let $G=(V,E)$ be a graph of order $n$ and let $\gamma _{R}(G)$
and $\partial (G)$ denote the Roman domination number and the differential
of $G,$ respectively. In this paper we prove that for any integer $k\geq 0$,
if $G$ is a graph of order $n\geq 6k+9$, minimum degree $\delta \geq 2,$
which does not contain any induced $\{C_{5},C_{8},\ldots ,C_{3k+2}\}$%
-cycles, then $\gamma _{R}(G)\leq \frac{(4k+8)n}{6k+11}$. This bound is an
improvement of the bounds given in [E.W. Chambers, B. Kinnersley, N. Prince,
and D.B. West, Extremal problems for Roman domination, SIAM J. Discrete
Math. 23 (2009) 1575--1586] when $k=0,$ {and [S. Bermudo, On the
differential and Roman domination number of a graph with minimum degree two,
Discrete Appl. Math. 232 (2017), 64--72] when }$k=1.$ Moreover, using the
Gallai-type result involving the Roman domination number and the
differential of graphs established by Bermudo et al. stating that $\gamma
_{R}(G)+\partial (G)=n$, we have $\partial (G)\geq \frac{(2k+3)n}{6k+11},$
thereby settling the conjecture of Bermudo posed in the second paper.

\noindent \textbf{Keyword:} Differential of a graph, Roman domination number.%
\newline
\textbf{MSC 2010}: 05C69
\end{abstract}



\section{Introduction}

In this paper, $G$ is a simple graph without isolated vertices, with vertex
set $V=V(G)$ and edge set $E=E(G)$. The \emph{order} $|V|$ of $G$ is denoted
by $n=n(G)$. For a vertex $v\in V$, the \emph{open neighborhood} of $v$ is
the set $N(v)=\{u\in V\mid uv\in E\},$ the \emph{closed neighborhood} of $v$
is the set $N[v]=N(v)\cup \{v\}$, and the \emph{degree} of $v$ is $\deg
_{G}(v)=|N(v)|$. 
Let $u$ and $v$ be two vertices in $G.$ A $uv$\textit{-path} is a path with
endvertices $u$ and $v$, and the \textit{distance} between $u$ and $v$ is
the length of a shortest $uv$-path. The \emph{diameter} of $G$, denoted by $%
\mathrm{diam}(G)$, is the maximum distance between vertices of $G$. We write
$P_{n}$ and $C_{n}$ for the \emph{path} and \emph{cycle} of order $n$,
respectively. {Let }$A$ and $B$ {are} two disjoint subgraphs {%
(not necessarily induced) }of a graph $G.$ If there is an edge $%
e $ having one endvertex in $A$ and the other one in $B,$ then $A+B+e$ will
denote the graph formed by $A$ and $B$ for which we add only the edge $e.$
We also denote by $G-A$ the subgraph of $G$ induced by $V(G)-V(A).$%

For a set $D,$ let $B(D)$ be the set of vertices in $V\setminus D$ that have
a neighbor in $D$. The \textit{differential} of a set $D$ is defined in~\cite%
{M} as $\partial (D)=|B(D)|-|D|$, and the maximum value of $\partial (D)$
for any subset $D$ of $V$ is the \textit{differential} of $G$, denoted $%
\partial (G)$. Differential of graphs has been studied extensively in
several papers, in particular \cite{B0,B3,B4,B5,B6,B1,B2,kkcs}. In 2017,
Bermudo \cite{B0} proved that for any graph $G$ with order $n\geq 15$,
minimum degree two and without any induced tailed $5$-cycle graph of seven
vertices or tailed $5$-cycle graph of seven vertices together with a
particular edge, it is satisfied $\partial (G)\geq \frac{5n}{17}$. Moreover,
he posed the following conjecture.

\begin{conj}[\protect\cite{B0}]
\label{conj}Let $G$ be a graph of order $n\ge 6k + 9$, minimum degree $%
\delta\ge 2$, which does not contain any induced $\{C_5,C_8,\ldots,C_{3k+2}%
\} $-cycles. Then $\partial(G)\ge \frac{(2k+3)n}{6k+11}.$
\end{conj}

A \emph{Roman dominating function} (RDF) on a graph $G$ is a function $%
f:V(G)\rightarrow \{0,1,2\}$ such that every vertex $u\in V(G)$ with $f(u)=0$
has a neighbor $v$ with $f(v)=2$. The \textit{weight} of an RDF $f$ is the
value $f(V(G))=\sum_{u\in V(G)}f(u),$ and the \emph{Roman domination number}
$\gamma _{R}(G)$ is the minimum weight of an RDF on $G$. The Roman
domination number of graphs was introduced in 2004 by Cockayne et al.\ in
\cite{CDHH} and is now well-studied in graph theory. The literature on Roman
domination and its variations has been surveyed and detailed in two book
chapters and three surveys \cite{Ch1, Ch2, Ch3, Ch4, Ch5}.

{In \cite{cham}, it has been shown that, if $G$ is a graph of
order $n\geq 9$ and minimum degree $\delta \geq 2$, then $\gamma _{R}(G)\leq
\frac{8n}{11}$. It was also shown in \cite{B0} that $\gamma _{R}(G)\leq
\frac{12n}{17}$ for any graph }$G$ {with order $n\geq 15$, minimum degree
two and without any induced tailed 5-cycle graph of seven vertices or tailed
5-cycle graph of seven vertices together with a particular edge.  }

In this paper, we improve the aforementioned known results by showing that
if $G$ is a graph satisfying the statement of Conjecture \ref{conj}, then $%
\gamma _{R}(G)\leq \frac{(4k+8)n}{6k+11}.$ U{sing the Gallai-type result
involving the differential and the Roman domination number of graphs
established by Bermudo, Fernau and Sigarreta \cite{B6} who proved that $%
\gamma _{R}(G)+\partial (G)=n$, our bound leads to }$\partial (G)\geq \frac{%
(2k+3)n}{6k+11}$ which settles Conjecture {\ref{conj}.  }

We close this section by recalling the exact values of the Roman domination
number of paths and cycles given in \cite{CDHH}, namely $\gamma
_{R}(P_{n})=\gamma _{R}(C_{n})=\lceil \frac{2n}{3}\rceil .$

\section{Some useful lemmas}

We gather in this section some results that will be useful to us thereafter.
For technical reasons, we will often consider three Roman dominating
functions $f_{1},f_{2}$, and $f_{3}$ on a graph $G$, where we use $%
\overrightarrow{f}$ to denote the $3$-tuple $(f_{1},f_{2},f_{3})$, and $%
\overrightarrow{f}(v)$ for $(f_{1}(v),f_{2}(v),f_{3}(v))$ for a vertex $v$.
A vertex $v$ is said to be $\overrightarrow{f}$\textit{-strong} if $%
f_{j}(v)=2$ for some $j\in \{1,2,3\}$. Moreover, the weight of $%
\overrightarrow{f}$ is $\omega (\overrightarrow{f})=\sum_{j=1}^{3}\omega
(f_{j})$. Clearly, $\omega (f_{j})\leq \omega (\overrightarrow{f})/3$ for
some $j\in \{1,2,3\}$. \ Also, if $H$ is an induced subgraph of $G$ and $f$
an RDF on $G$, then we denote the restriction of $f$ on $H$ by $f|_{V(H)}$
and let $f(V(H))=\omega (f,H).$

For integers $m$ and {$\ell$} such that $m\geq 3$ and $\ell \geq 1$, let $%
C_{m,\ell }$ be the graph obtained from a cycle {$C_{m}=x_{1}x_{2}\ldots
x_{m}x_{1}$} and a path $P=y_{1}y_{2}\ldots y_{\ell }$ by adding the edge $%
x_{1}y_{1},$ with {$y_{i}\notin V(C_{m})$} for all possible $i$. The graph $%
C_{m,\ell }$ will be called a \textit{tailed }$m$\textit{-cycle graph} of
order $m+\ell $. {We call an \textit{ear }of a cycle }$C$ in a graph $G,$ {to a
path }$P$ {in }$G-C$ {whose endvertices are adjacent to some vertices in }$%
C. $

\begin{lemma}
\label{ear1}Let $G$ be a graph, $u,v\in V(G)$ and $\overrightarrow{f}%
=(f_{1},f_{2},f_{3})$ be a 3-tuple of RDFs of $G$ such that $u$ and $v$ are $%
\overrightarrow{f}$-strong. If $H$ is a graph obtained from $G$ by adding a
path $Q=y_{1}\ldots y_{\ell }$ and the edges $uy_{1},vy_{\ell }$, then $%
\overrightarrow{f}$ can be extended to a 3-tuple of RDFs $\overrightarrow{g}$
of $H$ such that {$\omega (\overrightarrow{g},Q)\leq 2\ell$} and each vertex in $%
V(Q)-\{y_{1},y_{\ell }\}$ is $\overrightarrow{g}$-strong.
\end{lemma}

\noindent \textbf{Proof.} By assumption, $f_{i}(u)=2$ and $f_{j}(v)=2$ for
some $i,j\in \{1,2,3\}$. Let us consider the following two cases.

\smallskip \noindent \textbf{Case 1.} $i\neq j$.\newline
Assume, without loss of generality, that $i=1$ and $j=2$. Consider the
following situations.

\smallskip \textbf{Subcase 1.1.} $\ell \equiv 0\pmod 3$.\newline
Define the functions $g_{1},g_{2}$ and $g_{3}$ on $V(H)$ as follows: $%
g_{1}(z)=f_{1}(z)$ for all $z\in V(G)$, $g_{1}(y_{3i+3})=2$ for $0\leq i\leq
\frac{\ell }{3}-1$, and $g_{1}(z)=0$ otherwise; $g_{2}(z)=f_{2}(z)$ for all $%
z\in V(G)$, $g_{2}(y_{3i+1})=2$ for $0\leq i\leq \frac{\ell }{3}-1$, and $%
g_{2}(z)=0$ otherwise; and $g_{3}(z)=f_{3}(z)$ for all $z\in V(G)$, $%
g_{3}(y_{3i+2})=2$ for $0\leq i\leq \frac{\ell }{3}-1$, and $g_{3}(z)=0$
otherwise.

\smallskip \textbf{Subcase 1.2.} $\ell \equiv 1\pmod 3$.\newline
Define the functions $g_{1},g_{2}$ and $g_{3}$ on $V(H)$ as follows: $%
g_{1}(z)=f_{1}(z)$ for all $z\in V(G)$, $g_{1}(y_{3i+3})=2$ for $0\leq i\leq
\frac{\ell -4}{3}$, and $g_{1}(z)=0$ otherwise; $g_{2}(z)=f_{2}(z)$ for all $%
z\in V(G)$, $g_{2}(y_{3i+2})=2$ for $0\leq i\leq \frac{\ell -4}{3}$, and $%
g_{2}(z)=0$ otherwise; and $g_{3}(z)=f_{3}(z)$ for $z\in V(G)$, $%
g_{3}(y_{3i+1})=2$ for $0\leq i\leq \frac{\ell -1}{3}$, and $g_{3}(z)=0$
otherwise.

\smallskip \textbf{Subcase 1.3.} $\ell \equiv 2\pmod 3$.\newline
Define the functions $g_{1},g_{2}$ and $g_{3}$ on $V(H)$ as follows: $%
g_{1}(z)=f_{1}(z)$ for all $z\in V(G)$, $g_{1}(y_{3i+3})=2$ for $0\leq i\leq
\frac{\ell -5}{3}$, $g_{1}(y_{\ell })=1$ and $g_{1}(z)=0$ otherwise; $%
g_{2}(z)=f_{2}(z)$ for all $z\in V(G)$, $g_{2}(y_{3i+2})=2$ for $0\leq i\leq
\frac{\ell -5}{3}$, $g_{2}(y_{\ell -1})=1$ and $g_{2}(z)=0$ otherwise; and $%
g_{3}(z)=f_{3}(z)$ for all $z\in V(G)$, $g_{3}(y_{3i+1})=2$ for $0\leq i\leq
\frac{\ell -2}{3}$, and $g_{3}(z)=0$ otherwise.

\bigskip

In either subcase, $g_{1},g_{2},g_{3}$ are RDFs of $H$ and thus $%
g=(g_{1},g_{2},g_{3})$ is a 3-tuple of RDFs of $H.$ In addition, $\omega
(\overrightarrow{g},Q)\leq 2\ell$ and each vertex of $V(Q)-\{y_{1},y_{\ell }\}$ is $%
\overrightarrow{g}$-strong.

\smallskip \noindent \textbf{Case 2.} $i=j$.\newline
Assume, without loss of generality, that $i=j=1$. Consider again the
following situations.

\smallskip \textbf{Subcase 2.1.} $\ell \equiv 0\pmod 3$.\newline
Define the functions $g_{1},g_{2}$ and $g_{3}$ on $V(H)$ as follows: $%
g_{1}(z)=f_{1}(z)$ for all $z\in V(G)$, $g_{1}(y_{2})=1$, $g_{1}(y_{3i+4})=2$
for $0\leq i\leq \frac{\ell -6}{3}$, and $g_{1}(z)=0$ otherwise; $%
g_{2}(z)=f_{2}(z)$ for all $z\in V(G)$, $g_{2}(y_{3i+2})=2$ for $0\leq i\leq
\frac{\ell }{3}-1$, and $g_{2}(z)=0$ otherwise; $g_{3}(z)=f_{3}(z)$ for all $%
z\in V(G)$, $g_{3}(y_{1})=1$, $g_{3}(y_{3i+3})=2$ for $0\leq i\leq \frac{%
\ell }{3}-1$, and $g_{3}(z)=0$ otherwise.

\smallskip \textbf{Subcase 2.2.} $\ell \equiv 1\pmod 3$.\newline
Define the functions $g_{1},g_{2}$ and $g_{3}$ on $V(H)$ as follows: $%
g_{1}(z)=f_{1}(z)$ for $z\in V(G)$, $g_{1}(y_{2})=g_{1}(y_{\ell -1})=1$, $%
g_{1}(y_{3i+4})=2$ for $0\leq i\leq \frac{\ell -7}{3}$, and $g_{1}(z)=0$
otherwise; $g_{2}(z)=f_{2}(z)$ for all $z\in V(G)$, $g_{2}(y_{1})=1$, $%
g_{2}(y_{3i+3})=2$ for $0\leq i\leq \frac{\ell -4}{3}$, and $g_{2}(z)=0$
otherwise; $g_{3}(z)=f_{3}(z)$ for all $z\in V(G)$, $g_{3}(y_{\ell })=1$, $%
g_{3}(y_{3i+2})=2$ for $0\leq i\leq \frac{\ell -4}{3}$, and $g_{3}(z)=0$
otherwise.

\smallskip \textbf{Subcase 2.3.} $\ell \equiv 2\pmod 3$.\newline
Define the functions $g_{1},g_{2}$ and $g_{3}$ on $V(H)$ as follows: $%
g_{1}(z)=f_{1}(z)$ for all $z\in V(G)$, $g_{1}(y_{3i+3})=2$ for $0\leq i\leq
\frac{\ell -5}{3}$ and $g_{1}(z)=0$ otherwise; $g_{2}(z)=f_{2}(z)$ for all $%
z\in V(G)$, $g_{2}(y_{3i+1})=2$ for $0\leq i\leq \frac{\ell -2}{3}$, and $%
g_{2}(z)=0$ otherwise; $g_{3}(z)=f_{3}(z)$ for all $z\in V(G)$, $%
g_{3}(y_{3i+2})=2$ for $0\leq i\leq \frac{\ell -2}{3}$, and $g_{3}(z)=0$
otherwise.

\bigskip

In either subcase, $g_{1},g_{2},g_{3}$ are RDFs of $H$ and thus $%
g=(g_{1},g_{2},g_{3})$ is a 3-tuple of RDFs of $H.$ Moreover, $\omega
(\overrightarrow{g},Q)\leq 2\ell$ and each vertex of $V(Q)-\{y_{1},y_{\ell }\}$ is $%
\overrightarrow{g}$-strong.$\hfill \Box $

\begin{lemma}
\label{tailedcycle-3p+1}Let $G$ be a graph, $u\in V(G)$ and $\overrightarrow{%
f}=(f_{1},f_{2},f_{3})$ a 3-tuple of RDFs of $G$ such that $u$ is $%
\overrightarrow{f}$-strong.

\begin{itemize}
\item[1.] If $H$ is obtained from $G$ by adding a cycle $C_{3p+1}=x_{1}x_{2}%
\ldots x_{3p+1}x_{1}\;(p\geq 1)$ and the edge $ux_{1}$, then $%
\overrightarrow{f}$ can be extended to a 3-tuple $\overrightarrow{g}$ of
RDFs of $H$ such that {$\omega (\overrightarrow{g},C_{3p+1})\leq 2(3p+1)$} and each vertex
in $V(C_{3p+1})-\{x_{3p+1}\}$ is $\overrightarrow{g}$-strong.

\item[2.] If $H$ is obtained from $G$ by adding a tailed cycle $C_{3p+1,\ell
}\;(p\geq 1)$ and the edge $uy_{\ell }$, then $\overrightarrow{f}$ can be
extended to a 3-tuple $\overrightarrow{g}$ of RDFs of $H$ such that $\omega
(\overrightarrow{g},C_{3p+1,\ell })\leq 2(3p+1+\ell )$ and each vertex of $%
V(C_{3p+1,\ell })-\{x_{3p+1}\}$ is $\overrightarrow{g}$-strong.

\item[3.] If $H$ is obtained from $G$ by adding a cycle $C_{3p+2}=x_{1}x_{2}%
\ldots x_{3p+2}x_{1}\;(p\geq 1)$ and the edge $ux_{1}$, then $%
\overrightarrow{f}$ can be extended to a 3-tuple $\overrightarrow{g}$ of
RDFs of $H$ such that $\omega (\overrightarrow{g},C_{3p+2})\leq 2(3p+2)+1$ and each
vertex of $C_{3p+2}$, is $\overrightarrow{g}$-strong.

\item[4.] If $H$ is obtained from $G$ by adding a tailed cycle $C_{3p+2,\ell
}\;(p\geq 1)$ and the edge $uy_{\ell }$, then $\overrightarrow{f}$ can be
extended to a 3-tuple $\overrightarrow{g}$ of RDFs of $H$ such that $\omega
(\overrightarrow{g},C_{3p+2,\ell })\leq 2(3p+2+\ell )+1$ and each vertex of $%
C_{3p+2,\ell }$, is $\overrightarrow{g}$-strong.
\end{itemize}
\end{lemma}

\noindent \textbf{Proof.} Since $u$ is $\overrightarrow{f}$-strong, let us
assume, without loss of generality, that $f_{1}(u)=2$.

1) Define the functions $g_{1},g_{2}$ and $g_{3}$ on $V(H)$ as follows: $%
g_{1}(z)=f_{1}(z)$ for all $z\in V(G)$, $g_{1}(x_{3i+3})=2$ for $0\leq i\leq
p-1$, and $g_{1}(z)=0$ otherwise; $g_{2}(z)=f_{2}(z)$ for all $z\in V(G)$, $%
g_{2}(x_{3p})=1$, $g_{2}(x_{3i+1})=2$ for $0\leq i\leq p-1$, and $g_{2}(z)=0$
otherwise; $g_{3}(z)=f_{3}(z)$ for all $z\in V(G)$, $g_{3}(x_{3p+1})=1$, $%
g_{3}(x_{3i+2})=2$ for $0\leq i\leq p-1$, and $g_{3}(z)=0$ otherwise.

Clearly, $g_{1},g_{2},g_{3}$ are RDFs of $H$ and thus $g=(g_{1},g_{2},g_{3})$
is a 3-tuple of RDFs of $H.$ In addition, $\omega (\overrightarrow{g},C_{3p+1})\leq
2(3p+1)$ and each vertex of $V(C_{3p+1})-\{x_{3p+1}\}$ is $%
\overrightarrow{g}$-strong.

2) Consider the following cases.

\smallskip \noindent \textbf{Case 1.} $\ell \equiv 0\pmod 3$.\newline
Define the functions $g_{1},g_{2}$ and $g_{3}$ on $V(H)$ as follows: $%
g_{1}(z)=f_{1}(z)$ for all $z\in V(G)$, $g_{1}(x_{3i+3})=2$ for $0\leq i\leq
p-1$, $g_{1}(y_{3i+1})=2$ for $0\leq i\leq \frac{\ell }{3}-1$, and $%
g_{1}(z)=0$ otherwise; $g_{2}(z)=f_{2}(z)$ for all $z\in V(G)$, $%
g_{2}(x_{3p})=1$, $g_{2}(x_{3i+1})=2$ for $0\leq i\leq p-1$, $%
g_{2}(y_{3i+3})=2$ for $0\leq i\leq \frac{\ell }{3}-1$, and $g_{2}(z)=0$
otherwise; $g_{3}(z)=f_{3}(z)$ for all $z\in V(G)$, $g_{3}(x_{3p+1})=1$, $%
g_{3}(x_{3i+2})=2$ for $0\leq i\leq p-1$, $g_{3}(y_{3i+2})=2$ for $0\leq
i\leq \frac{\ell }{3}-1$, and $g_{3}(z)=0$ otherwise.

\smallskip \noindent \textbf{Case 2.} $\ell \equiv 1\pmod 3$.\newline
Define the functions $g_{1},g_{2}$ and $g_{3}$ on $V(H)$ as follows: $%
g_{1}(z)=f_{1}(z)$ for all $z\in V(G)$, $g_{1}(x_{3p+1})=1$, $%
g_{1}(x_{3i+2})=2$ for $0\leq i\leq p-1$, $g_{1}(y_{3i+2})=2$ for $0\leq
i\leq \frac{\ell -4}{3}$, and $g_{1}(z)=0$ otherwise; $g_{2}(z)=f_{2}(z)$
for all $z\in V(G)$, $g_{1}(x_{3p})=1$, $g_{2}(x_{3i+1})=2$ for $0\leq i\leq
p-1$, $g_{2}(y_{3i+3})=2$ for $0\leq i\leq \frac{\ell -4}{3}$, and $%
g_{2}(z)=0$ otherwise; $g_{3}(z)=f_{3}(z)$ for all $z\in V(G)$, $%
g_{3}(x_{3i+3})=2$ for $0\leq i\leq p-1$, $g_{3}(y_{3i+1})=2$ for $0\leq
i\leq \frac{\ell -1}{3}$, and $g_{3}(z)=0$ otherwise.

\smallskip \noindent \textbf{Case 3.} $\ell \equiv 2\pmod 3$.\newline
Define the functions $g_{1},g_{2}$ and $g_{3}$ on $V(H)$ as follows: $%
g_{1}(z)=f_{1}(z)$ for all $z\in V(G)$, $g_{1}(x_{3p})=1$, $%
g_{1}(x_{3i+1})=2 $ for $0\leq i\leq p-1$, $g_{1}(y_{3i+3})=2$ for $0\leq
i\leq \frac{\ell -5}{3}$, and $g_{1}(z)=0$ otherwise; $g_{2}(z)=f_{2}(z)$
for all $z\in V(G)$, $g_{2}(x_{3p+1})=1,$ $g_{2}(x_{3i+2})=2$ for $0\leq
i\leq p-1$, $g_{2}(y_{3i+2})=2$ for $0\leq i\leq \frac{\ell -2}{3}$, and $%
g_{2}(z)=0$ otherwise; $g_{3}(z)=f_{3}(z)$ for all $z\in V(G)$, $%
g_{3}(x_{3i+3})=2$ for $0\leq i\leq p-1$, $g_{3}(y_{3i+1})=2$ for $0\leq
i\leq \frac{\ell -2}{3}$, and $g_{3}(z)=0$ otherwise.

\bigskip

In either case, $g_{1},g_{2},g_{3}$ are RDFs of $H$ and thus $%
\overrightarrow{g}=(g_{1},g_{2},g_{3})$ is a 3-tuple of RDFs of $H.$ In
addition, $\omega (\overrightarrow{g},C_{3p+1,\ell })\leq 2(3p+1+\ell) $
and each vertex of $V(C_{3p+1,\ell })-\{x_{3p+1}\}$ is $\overrightarrow{g}$%
-strong.

The proofs of the remaining items are similar and therefore omitted. $\hfill
\Box $

\begin{lemma}
\label{cyclepath}

\begin{enumerate}
\item \emph{\ Let $C=v_{1}v_{2}\ldots v_{t}v_{1}$ be a cycle on $t\geq 4$
vertices with $t\equiv 1\pmod 3$. Then $C$ has a 3-tuple of RDFs $%
\overrightarrow{f}=(f_{1},f_{2},f_{3})$ such that $\omega (\overrightarrow{f}%
)\leq 2t+1$ and all vertices of $C$ but $v_{t}$ are $\overrightarrow{f}$%
-strong. }

\item \emph{Let $C=v_{1}v_{2}\ldots v_{t}v_{1}$ be a cycle on $t\geq 3$
vertices with $t\equiv 0\pmod 3$. Then $C$ has a 3-tuple of RDFs $%
\overrightarrow{f}=(f_{1},f_{2},f_{3})$ such that $\omega (\overrightarrow{f}%
)\leq 2t$ and all vertices of $C$ are $\overrightarrow{f}$-strong.}

\item \emph{Let $C_{m,\ell }$ be a tailed $m$-cycle with $m\equiv 1\pmod
3$ and $V(C_{m,\ell })=\{x_{1},x_{2},\ldots ,x_{m},y_{1},y_{2},\ldots
,y_{\ell }\},$ where the $x_{i}$'s induce a cycle {$C_{m}$ and the }$y_{i}$%
's induce a path $P_{\ell }$. Then $G$ has a 3-tuple of RDFs $%
\overrightarrow{f}=(f_{1},f_{2},f_{3})$ such that $\omega (\overrightarrow{f}%
)\leq 2(m+\ell )+1$ and all vertices of $C_{m,\ell }$ but $x_{m}$ are $%
\overrightarrow{f}$-strong.}
\end{enumerate}
\end{lemma}

\emph{\ }

\noindent \textbf{Proof. }1) Define the functions $f_{1},f_{2}$ and $f_{3}$
on $V(C)$ as follows: $f_{1}(v_{t-1})=1,f_{1}(v_{3i+1})=2$ for $0\leq i\leq
\frac{t-4}{3}$ and $f_{1}(x)=0$ otherwise; $f_{2}(v_{t})=1,$ $%
f_{2}(v_{3i+2})=2$ for $0\leq i\leq \frac{t-4}{3}$ and $f_{2}(x)=0$
otherwise; $f_{3}(v_{1})=1,f_{3}(v_{3i+3})=2$ for $0\leq i\leq \frac{t-4}{3}$
and $f_{3}(x)=0$ otherwise. Clearly, $f_{1},f_{2},f_{3}$ are RDFs of $C.$
Hence $\overrightarrow{f}=(f_{1},f_{2},f_{3})$ is a 3-tuple of RDFs of $C$,
with $\omega (\overrightarrow{f})\leq 2t+1$ and all vertices of $C$ except $%
v_{t}$ are $\overrightarrow{f}$-strong.

2) Define the functions $f_{1},f_{2}$ and $f_{3}$ on $V(C)$ as follows: $%
f_{1}(v_{3i+1})=2$ for $0\leq i\leq \frac{t-3}{3}$ and $f_{1}(x)=0$
otherwise; $f_{2}(v_{3i+2})=2$ for $0\leq i\leq \frac{t-3}{3}$ and $%
f_{2}(x)=0$ otherwise; $f_{3}(v_{3i+3})=2$ for $0\leq i\leq \frac{t-3}{3}$
and $f_{3}(x)=0$ otherwise. Clearly, $f_{1},f_{2},f_{3}$ are RDFs of $C$.
Hence $\overrightarrow{f}=(f_{1},f_{2},f_{3})$ is a 3-tuple of RDFs of $C$
with $\omega (\overrightarrow{f})\leq 2t$ and each vertex of $C$ are $%
\overrightarrow{f}$-strong.

3) Define the functions $f_{1},f_{2}$ and $f_{3}$ on $V(C_{m,\ell }).$ For
vertices on $C_{m}$ as follows: $f_{1}(x_{m-1})=1,f_{1}(x_{3i+1})=2$ for $%
0\leq i\leq \frac{m-4}{3}$ and $f_{1}(x)=0$ otherwise; $f_{2}(x_{m})=1,$ $%
f_{2}(x_{3i+2})=2$ for $0\leq i\leq \frac{m-4}{3}$ and $f_{2}(x)=0$
otherwise; $f_{3}(x_{3i+3})=2$ for $0\leq i\leq \frac{m-4}{3}$ and $%
f_{3}(x)=0$ otherwise. Moreover, the $f_{i}$'s are defined for the vertices
on $P_{\ell }$ according to $\ell $ as follows.

If $\ell \equiv 0\pmod 3$, then $f_{1}(y_{3i+3})=2$ for $0\leq i\leq \frac{%
\ell -3}{3}$ and $f_{1}(x)=0$ otherwise; $f_{2}(y_{3i+2})=2$ for $0\leq
i\leq \frac{\ell -3}{3}$ and $f_{2}(x)=0$ otherwise; $f_{3}(y_{\ell })=1$, $%
f_{3}(y_{3i+1})=2$ for $0\leq i\leq \frac{\ell -3}{3}$ and $f_{3}(x)=0$
otherwise.

If $\ell \equiv 1\pmod 3$, then $f_{1}(y_{3i+3})=2$ for $0\leq i\leq \frac{%
\ell -4}{3}$ and $f_{1}(x)=0$ otherwise; $f_{2}(y_{\ell })=1,$ $%
f_{2}(y_{3i+2})=2$ for $0\leq i\leq \frac{\ell -4}{3}$ and $f_{2}(x)=0$
otherwise; $f_{3}(y_{3i+1})=2$ for $0\leq i\leq \frac{\ell -1}{3}$ and $%
f_{3}(x)=0$ otherwise.

If $\ell \equiv 2\pmod 3$, then $f_1(y_{\ell })=1,$ $f_{1}(y_{3i+3})=2$ for $%
0\leq i\leq \frac{\ell -5}{3}$ and $f_{1}(x)=0$ otherwise; $%
f_{2}(y_{3i+2})=2 $ for $0\leq i\leq \frac{\ell -2}{3}$ and $f_{2}(x)=0$
otherwise; $f_{3}(y_{3i+1})=2$ for $0\leq i\leq \frac{\ell -2}{3}$ and $%
f_{3}(x)=0$ otherwise.

Clearly, in either case $f_{1},f_{2},f_{3}$ are RDFs of {$C_{m,\ell}$} and thus $%
\overrightarrow{f}=(f_{1},f_{2},f_{3})$ is a 3-tuple of RDFs of {$C_{m,\ell}.$} Also, $%
\omega (\overrightarrow{f})\leq 2(m+\ell )+1$ and all vertices of $C_{m,\ell
}$ but $x_{m}$ are $\overrightarrow{f}$-strong.$\hfill \Box $

\begin{lemma}
\label{MainLem}Let $C_{i}=x_{1}^{i}x_{2}^{i}\ldots x_{n_{i}}^{i}x_{1}^{i}\;$%
be a cycle of order $n_{i},$ for $i\in \{1,2,\ldots ,s\}.$

\begin{enumerate}
\item If $n_{1}\equiv 2\pmod 3$ and $G$ is a graph obtained from $C_{1}$ and
$C_{2}$ by identifying the vertices $x_{1}^{1}$ and $x_{1}^{2}$, then $G$
has a 3-tuple of RDFs $\overrightarrow{f}$ such that $\omega (%
\overrightarrow{f})\leq 2n(G)+1$, and all vertices of $V(G)-%
\{x_{2}^{2},x_{n_{2}}^{2}\}$ are $\overrightarrow{f}$-strong.

\item If $n_{1}\equiv 2\pmod 3$, $n_{2}\equiv 1\pmod 3$ and $G$ is obtained
from $C_{1}$ and $C_{2}$ by adding either the edge $x_{1}^{1}x_{1}^{2}$ or a
path $z_{1}\ldots z_{3k}\;(k\geq 1)$ and the edges $x_{1}^{1}z_{1},$ $%
x_{1}^{2}z_{3k}$, then $G$ has a 3-tuple of RDFs $\overrightarrow{f}%
=(f_1,f_2,f_3)$ such that $\omega (\overrightarrow{f})\leq 2n(G)+1$ and each
vertex of $G$ but $x_{n_{2}}^{2}$ is $\overrightarrow{f}$-strong.

\item If $n_{1}\equiv 2\pmod 3$, $n_{2}\equiv 1\pmod 3$ and $G$ is obtained
from $C_{1}$ and $C_{2}$ by adding for $k\geq 1$, a path $z_{1}\ldots
z_{3k+1}$ and the edges $x_{1}^{1}z_{1},$ $x_{1}^{2}z_{3k+1}$, then $G$ has
a 3-tuple of RDFs $\overrightarrow{f}$ such that $\omega (\overrightarrow{f}%
)\leq 2n(G)+1$ and all vertices of $G$ but $x_{n_{2}}^{2}$ are $%
\overrightarrow{f}$-strong.

\item If $n_{1}\equiv 2\pmod 3$, $n_{2}\equiv 1\pmod 3$ and $G$ is obtained
from $C_{1}$ and $C_{2}$ by adding for $k\geq 1$ a path $z_{1}\ldots
z_{3k+2}\;$and the edges $x_{1}^{1}z_{1},$ $x_{1}^{2}z_{3k+2}$, then $G$ has
a 3-tuple of RDFs $\overrightarrow{f}$ such that $\omega (\overrightarrow{f}%
)\leq 2n(G)+1$ and all vertices of $G$ but $x_{n_{2}}^{2}$ are $%
\overrightarrow{f}$-strong.

\item If $n_{i}\equiv 2\pmod 3$ for $i\in \{1,2\}$ and $G$ is obtained from $%
C_{1}$ and $C_{2}$ by adding either the edge $x_{1}^{1}x_{1}^{2}$ or a path $%
z_{1}\ldots z_{3k}\;(k\geq 1)$ and the edges $x_{1}^{1}z_{1},$ $%
x_{1}^{2}z_{3k}$, then $G$ has a 3-tuple of RDFs $\overrightarrow{f}$ such
that $\omega (\overrightarrow{f})\leq 2n(G)+2$ and all vertices of $G$ are $%
\overrightarrow{f}$-strong.

\item If $n_{i}\equiv 2\pmod 3$ for $i\in \{1,2\}$ and $G$ is obtained from%
\textbf{\ }$C_{1}$ and $C_{2}$ by adding a path $z_{1}\ldots
z_{3k+1}\;(k\geq 1)$ and the edges $x_{1}^{1}z_{1},$ $x_{1}^{2}z_{3k+1}$,
then $G$ has a 3-tuple of RDFs $\overrightarrow{f}$ such that $\omega (%
\overrightarrow{f})\leq 2n(G)+2$ and all vertices of $G$ are $%
\overrightarrow{f}$-strong.

\item If $n_{i}\equiv 2\pmod 3$ for $i\in \{1,2\}$ and $G$ is obtained from $%
C_{1}$ and $C_{2}$ by adding a path $z_{1}\ldots z_{3k+2}\;(k\geq 1)$ and
the edges $x_{1}^{1}z_{1},$ $x_{1}^{2}z_{3k+2}$, then $G$ has a 3-tuple of
RDFs $\overrightarrow{f}$ such that $\omega (\overrightarrow{f})\leq 2n(G)+2$
and all vertices of $G$ are $\overrightarrow{f}$-strong.

\item If $n_{1}\equiv 2\pmod 3$, $n_{2}\equiv 0\pmod 3$ and $G$ is obtained
from $C_{1}$ and $C_{2}$ by adding either the edge $x_{1}^{1}x_{1}^{2}$ or a
path $z_{1}\ldots z_{3k}\;(k\geq 1)$ and the edges $x_{1}^{1}z_{1},$ $%
x_{1}^{2}z_{3k}$, then $G$ has a 3-tuple of RDFs $\overrightarrow{f}$ such
that $\omega (\overrightarrow{f})\leq 2n(G)+1$ and all vertices of $G$ are $%
\overrightarrow{f}$-strong.

\item If $n_{1}\equiv 2\pmod 3$, $n_{2}\equiv 0\pmod 3$ and $G$ is obtained
from $C_{1}$ and $C_{2}$ by adding a path $z_{1}\ldots z_{3k+1}\;(k\geq 1)$
and the edges $x_{1}^{1}z_{1},$ $x_{1}^{2}z_{3k+1}$, then $G$ has a 3-tuple
of RDFs $\overrightarrow{f}$ such that $\omega (\overrightarrow{f})\leq
2n(G)+1$ and all vertices of $G$ are $\overrightarrow{f}$-strong.

\item If $n_{1}\equiv 2\pmod 3$, $n_{2}\equiv 0\pmod 3$ and $G$ is obtained
from $C_{1}$ and $C_{2}$ by adding a path $z_{1}\ldots z_{3k+2}\;(k\geq 1)$
and the edges $x_{1}^{1}z_{1},$ $x_{1}^{2}z_{3k+2}$, then $G$ has a 3-tuple
of RDFs $\overrightarrow{f}$ such that $\omega (\overrightarrow{f})\leq
2n(G)+1$ and all vertices of $G$ are $\overrightarrow{f}$-strong.

\item {If $s\geq 3$, $n_{i}\equiv 2\pmod 3$ for each $i$ and $G$ is obtained
from $C_{1},\ldots ,C_{s}$ by adding a new vertex $x$ and the edges $%
xx_{1}^{1},\ldots ,xx_{1}^{s}$, then $G$ has a 3-tuple of RDFs $%
\overrightarrow{f}$ such that $\omega (f,G)\leq 2n(G)-s+4$ and the vertex $x$
is $\overrightarrow{f}$-strong.}
\end{enumerate}
\end{lemma}

\noindent\textbf{Proof.}

\begin{enumerate}
\item Define the functions $f_{1},f_{2}$ and $f_{3}$ on $V(G)$ as follows$.$
For vertices on $C_{1}$: $f_{1}(x_{3i+1}^{1})=f_{2}(x_{3i+2}^{1})=2$ for
each $0\leq i\leq \frac{n_{1}-2}{3}$ and $f_{i}(x)=0$ otherwise, {%
for $i=1,2$, and $f_{3}(x_{3i+3}^{1})=2$ for each $0\leq i\leq
\frac{n_{1}-5}{3}$, $f_{3}(x_{1}^{1})=2$, and $f_{3}(x)=0$ otherwise.} Now
for vertices on $C_{2}$ but $x_{1}^{2}$:

If $n_{2}\equiv 0\pmod 3$, then let $f_{1}(x_{3i+4}^{2})=2$ for each $0\leq
i\leq \frac{n_{2}-6}{3}$ and $f_{1}(x)=0$ otherwise; $f_{2}(x_{3i+3}^{2})=2$
for each $0\leq i\leq \frac{n_{2}-3}{3}$ and $f_{2}(x)=0$ otherwise; $%
f_{3}(x_{3i+2}^{2})=2$ for each $0\leq i\leq \frac{n_{2}-3}{3}$ and $%
f_{3}(x)=0$ otherwise.

If $n_{2}\equiv 1\pmod 3$, then let $f_{1}(x_{3i+4}^{2})=2$ for each $0\leq
i\leq \frac{n_{2}-7}{3}$, $f_{1}(x_{n_{2}-1}^{2})=1$ and $f_{1}(x)=0$
otherwise; $f_{2}(x_{3i+3}^{2})=2$ for each $0\leq i\leq \frac{n_{2}-4}{3}$
and $f_{2}(x)=0$ otherwise; {$f_{3}(x_{n_{2}}^{2})=1$,} $%
f_{3}(x_{3i+2}^{2})=2$ for each $0\leq i\leq \frac{n_{2}-4}{3}$ and $%
f_{3}(x)=0$ otherwise.

If $n_{2}\equiv 2\pmod 3$, then let $f_{1}(x_{3i+4}^{2})=2$ for each $0\leq
i\leq \frac{n_{2}-5}{3}$ and $f_{1}(x)=0$ otherwise; $f_{2}(x_{3i+3}^{2})=2$
for each $0\leq i\leq \frac{n_{2}-5}{3}$, $f_{2}(x_{n_{2}}^{2})=1$ and $%
f_{2}(x)=0$ otherwise;
$f_{3}(x_{3}^{2})=f_{3}(x_{n_{2}-1}^{2})=1$, $f_{3}(x_{3i+5}^{2})=2$ for each
$0\leq i\leq \frac{n_{2}-8}{3}$ and $f_{3}(x)=0$ otherwise.

In either case, $f_{1},f_{2}$ and $f_{3}$ are RDFs of $G.$ Hence $\omega (%
\overrightarrow{f})\leq 2n(G)+1$, and all vertices of $V(G)-%
\{x_{2}^{2},x_{n_{2}}^{2}\}$ are $\overrightarrow{f}$-strong. \textbf{\ }

\item Define the functions $f_{1},f_{2}$ and $f_{3}$ on $V(G)$ as follows: $%
f_{1}(x_{3i+1}^{1})=2$ for each $0\leq i\leq \frac{n_{1}-2}{3}$, $%
f_{1}(z_{3i+3})=2$ for $0\leq i\leq k-1$, $f_{1}(x_{3i+3}^{2})=2$ for each $%
0\leq i\leq \frac{n_{2}-4}{3}$ and $f_{1}(x)=0$ otherwise; $%
f_{2}(x_{3i+2}^{1})=2$ for each $0\leq i\leq \frac{n_{1}-2}{3}$, $%
f_{2}(z_{3i+2})=2$ for $0\leq i\leq k-1$, $f_{2}(x_{3i+2}^{2})=2$ for each $%
0\leq i\leq \frac{n_{2}-4}{3}$, $f_{2}(x_{n_{2}}^{2})=1$ and $f_{2}(x)=0$
otherwise; $f_{3}(x_{n_{1}}^{1})=1$, $f_{3}(x_{3i+3}^{1})=2$ for each $0\leq
i\leq \frac{n_{1}-5}{3}$, $f_{3}(z_{3i+1})=2$ for $0\leq i\leq k-1$, $%
f_{3}(x_{3i+1}^{2})=2$ for each $0\leq i\leq \frac{n_{2}-4}{3}$, $%
f_{3}(x_{n_{2}-1}^{2})=2$ and $f_{3}(x)=0$ otherwise.

Clearly $f_{1},f_{2}$ and $f_{3}$ are RDFs of $G,$ and thus $\overrightarrow{%
f}=(f_{1},f_{2},f_{3})$ is a 3-tuple of RDFs of $G.$ Moreover, $\omega (%
\overrightarrow{f})\leq 2n(G)+1$ and each vertex of $G$ but $x_{n_{2}}^{2}$
is $\overrightarrow{f}$-strong.

\item Define the functions $f_{1},f_{2}$ and $f_{3}$ on $V(G)$ as follows: $%
f_{1}(x_{3i+1}^{1})=2$ for each $0\leq i\leq \frac{n_{1}-2}{3}$, $%
f_{1}(z_{3i+3})=2$ for $0\leq i\leq k-1$, $f_{1}(x_{3i+2}^{2})=2$ for each $%
0\leq i\leq \frac{n_{2}-4}{3}$, $f_{1}(x_{n_{2}}^{2})=1$ and $f_{1}(x)=0$
otherwise; $f_{2}(x_{3i+2}^{1})=2$ for each $0\leq i\leq \frac{n_{1}-2}{3}$,
$f_{2}(z_{3i+2})=2$ for $0\leq i\leq k-1$, $f_{2}(x_{3i+1}^{2})=2$ for each $%
0\leq i\leq \frac{n_{2}-4}{3}$, $f_{2}(x_{n_{2}-1}^{2})=1$ and $f_{2}(x)=0$
otherwise; $f_{3}(x_{n_{1}}^{1})=1$, $f_{3}(x_{3i+3}^{1})=2$ for each $0\leq
i\leq \frac{n_{1}-5}{3}$, $f_{3}(z_{3i+1})=2$ for $0\leq i\leq k$, $%
f_{3}(x_{3i+3}^{2})=2$ for each $0\leq i\leq \frac{n_{2}-4}{3}$ and $%
f_{3}(x)=0$ otherwise.

Clearly $f_{1},f_{2}$ and $f_{3}$ are RDFs of $G$ and thus $\overrightarrow{f%
}=(f_{1},f_{2},f_{3})$ is a 3-tuple of RDFs of $G.$ Also, $\omega (%
\overrightarrow{f})\leq 2n(G)+1$ and all vertices of $G$ but $x_{n_{2}}^{2}$
are $\overrightarrow{f}$-strong.

\item Define the functions $f_{1},f_{2}$ and $f_{3}$ on $V(G)$ as follows: $%
f_{1}(x_{3i+1}^{1})=2$ for each $0\leq i\leq \frac{n_{1}-2}{3}$, $%
f_{1}(z_{3i+3})=2$ for $0\leq i\leq k-1$, $f_{1}(x_{3i+1}^{2})=2$ for each $%
0\leq i\leq \frac{n_{2}-4}{3}$, $f_{3}(x_{n_{2}-1}^{2})=1$ and $f_{1}(x)=0$
otherwise; $f_{2}(x_{3i+2}^{1})=2$ for each $0\leq i\leq \frac{n_{1}-2}{3}$,
$f_{2}(z_{3i+2})=2$ for $0\leq i\leq k$, $f_{2}(x_{3i+3}^{2})=2$ for each $%
0\leq i\leq \frac{n_{2}-4}{3}$ and $f_{2}(x)=0$ otherwise; $%
f_{3}(x_{n_{1}}^{1})=1$, $f_{3}(x_{3i+3}^{1})=2$ for each $0\leq i\leq \frac{%
n_{1}-5}{3}$, $f_{3}(z_{3i+1})=2$ for $0\leq i\leq k$, $%
f_{3}(x_{3i+2}^{2})=2 $ for each $0\leq i\leq \frac{n_{2}-4}{3}$, $%
f_{3}(x_{n_{2}}^{2})=1$ and $f_{3}(x)=0$ otherwise.

Clearly $f_{1},f_{2}$ and $f_{3}$ are RDFs of $G$ and thus $\overrightarrow{f%
}=(f_{1},f_{2},f_{3})$ is a 3-tuple of RDFs of $G.$ Also, $\omega (%
\overrightarrow{f})\leq 2n(G)+1$ and each vertex of $G$ but $x_{n_{2}}^{2}$
is $\overrightarrow{f}$-strong.

\item Define the functions $f_{1},f_{2}$ and $f_{3}$ on $V(G)$ as follows: $%
f_{1}(x_{3i+1}^{1})=2$ for each $0\leq i\leq \frac{n_{1}-2}{3}$, $%
f_{1}(z_{3i+3})=2$ for $0\leq i\leq k-1$, $f_{1}(x_{3i+3}^{2})=2$ for each $%
0\leq i\leq \frac{n_{2}-5}{3}$, $f_{1}(x_{n_{2}}^{2})=1$ and $f_{1}(x)=0$
otherwise; $f_{2}(x_{3i+2}^{1})=2$ for each $0\leq i\leq \frac{n_{1}-2}{3}$,
$f_{2}(z_{3i+2})=2$ for $0\leq i\leq k-1$, $f_{2}(x_{3i+2}^{2})=2$ for each $%
0\leq i\leq \frac{n_{2}-2}{3}$ and $f_{2}(x)=0$ otherwise; $%
f_{3}(x_{n_{1}}^{1})=1$, $f_{3}(x_{3i+3}^{1})=2$ for each $0\leq i\leq \frac{%
n_{1}-5}{3}$, $f_{3}(z_{3i+1})=2$ for $0\leq i\leq k-1$, $%
f_{3}(x_{3i+1}^{2})=2$ for each $0\leq i\leq \frac{n_{2}-2}{3}$ and $%
f_{3}(x)=0$ otherwise.

Clearly $f_{1},f_{2}$ and $f_{3}$ are RDFs of $G$ and thus $\overrightarrow{f%
}=(f_{1},f_{2},f_{3})$ is a 3-tuple of RDFs of $G.$ Further, $\omega (%
\overrightarrow{f})\leq 2n(G)+2$ and all vertices of $G$ are $%
\overrightarrow{f}$-strong.

\item The proof is similar to that of item (5).

\item The proof is similar to that of item (5).

\item Define the functions $f_{1},f_{2}$ and $f_{3}$ on $V(G)$ as follows: $%
f_{1}(x_{3i+1}^{1})=2$ for each $0\leq i\leq \frac{n_{1}-2}{3}$, $%
f_{1}(z_{3i+3})=2$ for $0\leq i\leq k-1$, $f_{1}(x_{3i+3}^{2})=2$ for each $%
0\leq i\leq \frac{n_{2}-3}{3}$ and $f_{1}(x)=0$ otherwise; $%
f_{2}(x_{3i+2}^{1})=2$ for each $0\leq i\leq \frac{n_{1}-2}{3}$, $%
f_{2}(z_{3i+2})=2$ for $0\leq i\leq k-1$, $f_{2}(x_{3i+2}^{2})=2$ for each $%
0\leq i\leq \frac{n_{2}-3}{3}$ and $f_{2}(x)=0$ otherwise; $%
f_{3}(x_{n_{1}}^{1})=1$, $f_{3}(x_{3i+3}^{1})=2$ for each $0\leq i\leq \frac{%
n_{1}-5}{3}$, $f_{3}(z_{3i+1})=2$ for $0\leq i\leq k-1$, $%
f_{3}(x_{3i+1}^{2})=2$ for each $0\leq i\leq \frac{n_{2}-3}{3}$ and $%
f_{3}(x)=0$ otherwise.

Then $f_{1},f_{2}$ and $f_{3}$ are RDFs of $G$ and thus $\overrightarrow{f}%
=(f_{1},f_{2},f_{3})$ is a 3-tuple of RDFs of $G$ with the desired property.

\item The proof is similar to that of item (8).

\item The proof is similar to that of item (8).

\item {Define the function $f_{1}$ by $f_{1}(x)=2$, $f_{1}(x_{n_{j}}^{j})=1$
for each $1\leq j\leq s$, $f_{1}(x_{i}^{j})=2$ for each $j$ and each $%
i\equiv 0\pmod 3$, and $f_{1}(y)=0$ otherwise, and set $\overrightarrow{f}%
=(f_{1},f_{1},f_{1})$. }Clearly $f_{1}$ is an RDF of $G$ and thus $%
\overrightarrow{f}=(f_{1},f_{2},f_{3})$ is a 3-tuple of RDFs of {$G$ such
that $\omega (\overrightarrow{f})\leq 2n(G)-s+4$ and the vertex $x$ is $%
\overrightarrow{f}$-strong as desired.} $\hfill \Box $
\end{enumerate}

\begin{lemma}
\label{ear2}Let $H$ be a graph obtained from a cycle $C_{3p+2}=x_{1}x_{2}%
\ldots x_{3p+2}x_{1}$ and a path $Q=y_{1}\ldots y_{\ell }$ where $\ell
\equiv 1\;\mathrm{or}\;2\pmod 3$ by adding the edge $y_{1}x_{1}$ and joining
$y_{\ell }$ to some vertices in $V(C_{3p+2})-\{x_{1}\}$ with the condition
that:

\begin{description}
\item[(a)] if $\ell\equiv 1\pmod 3$ and $y_{\ell}x_{j}\in E(H)$, then $%
j\not\equiv 2\pmod 3$,

\item[(b)] if $\ell\equiv 2\pmod 3$ and $y_{\ell}x_{j}\in E(H)$, then {%
$j\equiv 2\pmod 3$.}
\end{description}

Then there exists a 3-tuple $\overrightarrow{g}=(g_{1},g_{2},g_{3})$ of RDF
of $H$ such that $\omega (\overrightarrow{g},H)\leq 2n(H)+1$ and each vertex
of $H$ but $y_{1},y_{\ell }$ is $\overrightarrow{g}$-strong.
\end{lemma}

\noindent \textbf{Proof.} First let $\ell \equiv 1\pmod 3$ and $y_{\ell
}x_{j}\in E(H).$ Define the functions $g_{1},g_{2}$ and $g_{3}$ on $V(H)$ as
follows, depending on whether $j\equiv 0\pmod 3$ or $j\equiv 1\pmod 3.$

If $j\equiv 0\pmod 3,$ then let $g_{1}(x_{3i+1})=2$ for $0\leq i\leq p$, $%
g_{1}(y_{3i+3})=2$ for $0\leq i\leq \frac{\ell -4}{3}$, and $g_{1}(z)=0$
otherwise; $g_{2}(x_{3i+2})=2$ for $0\leq i\leq p$, $g_{2}(y_{\ell })=1$, $%
g_{2}(y_{3i+2})=2$ for $0\leq i\leq \frac{\ell -4}{3}$, and $g_{2}(z)=0$
otherwise; $g_{3}(x_{1})=2,$ $g_{3}(x_{3i+3})=2$ for $0\leq i\leq p-1$, $%
g_{3}(y_{3i+2})=2$ for $0\leq i\leq \frac{\ell -4}{3}$, and $g_{1}(z)=0$
otherwise.

If $j\equiv 1\pmod 3$, then let $g_{1}(x_{3i+1})=2$ for $0\leq i\leq p$, $%
g_{1}(y_{3i+3})=2$ for $0\leq i\leq \frac{\ell -4}{3}$, and $g_{1}(z)=0$
otherwise; $g_{2}(x_{3i+2})=2$ for $0\leq i\leq p$, $g_{2}(y_{\ell })=1$, $%
g_{2}(y_{3i+2})=2$ for $0\leq i\leq \frac{\ell -4}{3}$, and $g_{2}(z)=0$
otherwise; $g_{3}(x_{1})=2,$ $g_{3}(x_{3i+3})=2$ for $0\leq i\leq p-1$, $%
g_{3}(y_{3i+3})=2$ for $0\leq i\leq \frac{\ell -4}{3}$, and $g_{1}(z)=0$
otherwise.

Second, let $\ell \equiv 2\pmod 3$ and $y_{\ell }x_{j}\in E(H),$ where {%
$j\equiv 2\pmod 3$.} Define the functions $g_{1},g_{2}$ and $%
g_{3}$ on $V(H)$ as follows: $g_{1}(x_{p})=1,$ $g_{1}(x_{3i+3})=2$ for $%
0\leq i\leq p-1$, $g_{1}(y_{3i+1})=2$ for $0\leq i\leq \frac{\ell -2}{3}$,
and $g_{1}(z)=0$ otherwise; $g_{2}(x_{3i+2})=2$ for $0\leq i\leq p$, $%
g_{2}(y_{3i+2})=2$ for $0\leq i\leq \frac{\ell -5}{3}$, $g_2(y_{\ell-1})=1$,
and $g_{2}(z)=0$ otherwise; $g_{3}(x_{3i+1})=2$ for $0\leq i\leq p$, $%
g_{3}(y_{3i+3})=2$ for $0\leq i\leq \frac{\ell -5}{3}$, $g_3(y_{\ell})=1$,
and $g_{1}(z)=0$ otherwise.

Clearly, $g_{1},g_{2},g_{3}$ are RDFs of $H$ and thus $\overrightarrow{g}%
=(g_{1},g_{2},g_{3})$ is a 3-tuple of RDFs of $H$ with the desired property.
$\hfill \Box $

\section{Partial answer to Conjecture \protect\ref{conj}}

In this section, we give a positive answer to Conjecture \ref{conj} for some
particular graphs. We start with the following two lemmas.

\begin{lemma}
\label{induced}\emph{Let }$k\geq 1$ \emph{be an integer and let $G$ be a
connected graph with $\delta \geq 2$, which does not contain {\ neither }any
induced $\{C_{5},C_{8},\ldots ,C_{3k+2}\}$-cycles {\ nor }any cycle of
length $\equiv 0\pmod 3$. Let $C$ be a cycle of $G$ with length $\ell (C)$ $%
\equiv 2\pmod 3$. Then }

\begin{enumerate}
\item \emph{if $C$ is induced in $G$, then $\ell (C)\geq 3k+5$, and }

\item \emph{if $C$ is not induced in $G$, then $\ell (C)\geq 6k+8$. }
\end{enumerate}
\end{lemma}

\noindent \textbf{Proof.} Item (1) is immediate since $G$ does not contain
any induced $\{C_{5},C_{8},\ldots ,C_{3k+2}\}$-cycles and $\ell (C)\equiv 2%
\pmod 3$. To prove item (2), let $C=v_{1}v_{2}\ldots v_{3p+2}v_{1}$ be a
cycle which is not induced in $G$. Hence $C$ has a chord, say without loss
of generality, $v_{1}v_{i}\in E(G)$. Consider the two paths $%
P=v_{i+1}v_{i+2}\ldots v_{3p+2}v_{1}$ and $Q=v_{2}v_{3}\ldots v_{i}$. Let $%
n(P)$ and $n(Q)$ denote the order of $P$ and $Q,$ respectively. Clearly $%
n(P)+n(Q)=3p+2$. Now, if $n(P)\equiv 0\pmod 3$, then $v_{1}v_{2}v_{3}\ldots
v_{i}v_{1}$ is a cycle of length $\equiv 0\pmod 3,$ contradicting the fact
that $G$ has no cycle of such length. Hence $n(P)\not\equiv 0\pmod 3$, and
likewise $n(Q)\not\equiv 0\pmod 3$. Moreover, since $n(P)+n(Q)=3p+2,$ we
deduce that $n(P)\not\equiv 2\pmod 3$ and $n(Q)\not\equiv 2\pmod 3$. Hence $%
n(P)\equiv 1\pmod 3$ and $n(Q)\equiv 1\pmod 3$. Consider the cycles $%
C_{1}=v_{1}v_{2}v_{3}\ldots v_{i}v_{1}$ and $C_{2}=v_{1}v_{i}v_{i+1}\ldots
v_{3p+2}v_{1}$. Then $\ell (C_{1})\equiv 2\pmod 3$ and $\ell (C_{2})\equiv 2%
\pmod 3$. If $C_{1}$ and $C_{2}$ are induced in $G$, then by item (1) we
have $i\geq 3k+5$ and $3p+4-i\geq 3k+5$ and thus $\ell (C)=3p+2\geq
3k+5+i-2\geq 3k+5+3k+5-2=6k+8$. Hence we assume that $C_{1}$ is not induced
in $G$. By repeating the above process we can see that the subgraph $%
G[V(C_{1})]$ has an induced cycle of length $\equiv 2\pmod 3$ and so $%
|V(C_{1})|=i\geq 3k+5$. If $C_{2}$ is an induced cycle, then by item (1) we
have $3p+4-i\geq 3k+5$ and so $\ell (C)=3p+2\geq 3k+5+i-2\geq
3k+5+3k+5-2=6k+8$. Now if $C_{2}$ is not an induced cycle, then a similar
argument as above shows that $G[V(C_{2})]$ has an induced cycle of length $%
\equiv 2\pmod 3$ yielding also $\ell (C)\geq 6k+8$.$\hfill \Box $

\begin{lemma}
\label{12}Let $G$ be a connected\textbf{\ }graph with minimum degree $\delta
\geq 2$ and let $G_{1}$ and $G_{2}$ be two non-null subgraphs of $G$ such
that $V(G)=V(G_{1})\cup V(G_{2})$. Then one of the following holds:

\begin{enumerate}
\item $G_{1}$ has a path $P=v_{1}\ldots v_{t}$ such that both $v_{1}$ and $%
v_{t}$ have neighbors in $G_{2}$ and $N_{G}(v_{1})\cup N_{G}(v_{t})\subseteq
V(G_{2})\cup V(P)$.

\item $G_{1}$ has a cycle $C=v_{1}v_{2}\ldots v_{t}v_{1}$ such that $v_{1}$
has neighbors in $G_{2}$ and $N_{G}(v_{t})\subseteq V(G_{2})\cup V(C)$.

\item $G_{1}$ contains a tailed $m$-cycle, say $C_{m,\ell }$, such that $%
y_{\ell }$ is adjacent to some vertex in $G_{2}$ and $N_{G}(x_{2})\cup
N_{G}(x_{m})\subseteq V(G_{2})\cup V(C_{m,\ell })$.
\end{enumerate}
\end{lemma}

\noindent \textbf{Proof.} Let $\mathcal{P}$ be the family of all longest
paths (not necessarily induced) in $G_{1}$ such that at least one of their
end-points has a neighbor in $G_{2}$ and let $Q=\{v\in V(G)\mid \mathrm{%
there\;is\;a\;path}$ $v_{1},\ldots ,v_{t}(=v)\in \mathcal{P}\;\mathrm{%
such\;that}\;v_{1}\;\mathrm{has}$ $\mathrm{a}$ $\mathrm{neighbor}$ $\mathrm{%
in}\;G_{2}\}$. Choose a vertex $v\in Q$ such that the length of its
corresponding path $P=v_{1},\ldots ,v_{t}(=v)\in \mathcal{P}$ is as long as
possible.

First let $v$ be adjacent to some vertex in $G_{2}$. By the definition of
set $Q,$ we have $v_{1}\in Q$, and from the choice of $v$ we deduce that $%
N_{G}(v_{1})\cup N_{G}(v_{t})\subseteq V(G_{2})\cup V(P).$ Hence item (1)
holds. Suppose now $v$ has no neighbor in $G_{2}$. It follows from the
choice of $v$ and the fact $\delta \geq 2$ that $v$ has at least two
neighbors in $V(P)$. Let $j$ be the smallest index such that $vv_{j}$ $\in
E(G).$ Now, if $j=1,$ that is $v$ is adjacent to $v_{1}$, then clearly $%
N_{G}(v_{t})\subseteq V(P)$ and thus item (2) holds. Hence we assume that $%
j\neq 1.$ Then $v_{1}\ldots v_{j}...v_{t-1}vv_{j}$\textbf{\ }is a tailed
cycle contained in $G_{1}.$ Observe that the path with endvertices $v_{1}$
and $v_{j+1}$ starting from $v_{1}$ to $v_{j}$ and then passing through $%
v_{t}$ to $v_{j+1}$ is also a longest path with same length as $P.$ Since $%
v=v_{t}$ has no neighbor in $G_{2},$ we may assume by analogy that $v_{j+1}$
has no neighbor in $G_{2}$ and thus all its neighbors are on\textbf{\ }$P$
which forms a tailed cycle and thus item (3) holds$.$\textbf{\ }$\hfill \Box
$

\begin{theorem}
\label{Th1}\emph{\emph{Let} $k\geq 1$ be an integer and let $G$ be a
connected graph of order $n\geq 6k+9$ and minimum degree at least 2 such
that $G$ has no cycle with length $\equiv 0\;\mathrm{or}\;2\pmod 3$. Then $%
\gamma _{R}(G)\leq \frac{(4k+8)n}{6k+11}$.}
\end{theorem}

\noindent \textbf{Proof.} Let $Q=z_{1}z_{2}\ldots z_{r}$ be a longest path
in $G$. If $V(G)=V(Q)$, then we have $\gamma _{R}(G)\leq \frac{2n+1}{3}<%
\frac{(4k+8)n}{6k+11}$. Hence, we assume that $V(Q)\subsetneqq V(G)$. By the
choice of $Q$ we have $N_{G}(z_{1})\cup N_{G}(z_{r})\subseteq V(Q)$. Since $%
\delta (G)\geq 2$, $z_{1}$ is adjacent to some $z_{j}$ with $j\equiv 1\pmod 3%
,$ because $G$ has no cycle with length $\equiv 0\;\mathrm{or}\;2\pmod 3$.
Let $G_{2}^{0}$ be the graph obtained from the path $Q$ to which we add the
edge $z_{1}z_{j}$ and let $G_{1}^{0}$ be the graph induced by $%
V(G)-V(G_{2}^{0}).$ Observe that $G_{2}^{0}$ is a tailed $j$-cycle $%
C_{j,r-j}.$ By Lemma \ref{cyclepath}-(3), $G_{2}^{0}$ has a 3-tuple of RDFs $%
\overrightarrow{f}=(f_{1},f_{2},f_{3})$ such that $\omega (\overrightarrow{f}%
)\leq 2r+1$ and all vertices of $C_{j,r-j}$ but $z_{1}$ are $\overrightarrow{%
f}$-strong. According to Lemma \ref{12}, we consider the following three
possibilities.

\begin{description}
\item[(a)] $G_{1}^{0}$ has a path $P=v_{1}\ldots v_{t}$ such that $%
v_{1},v_{t}$ are adjacent to some vertices in $V(G_{2}^{0})-\{z_{1},z_{t}\}$%
, say $u,v$ (possibly $u=v$) and $N_{G}(v_{1})\cup N_{G}(v_{t})\subseteq
V(G_{2}^{0})\cup V(P)$.\newline
Let $G_{2}^{1}$ be the graph obtained from $G_{2}^{0}$ and the path $P$ by
adding the edges $v_{1}u,v_{t}v$. By Lemma \ref{ear1}, $\overrightarrow{f}$
can be extended to a 3-triple of RDFs $\overrightarrow{g}$ of $G_{2}^{1}$
such that $\omega (\overrightarrow{g},P)\leq 2n(P)$ and each vertex in $%
V(P)-\{v_{1},v_{t}\}$ is $\overrightarrow{g}$-strong. Note that $\omega (%
\overrightarrow{g})=\omega (\overrightarrow{f})+\omega (\overrightarrow{g}%
,P)\leq 2n(G_{2}^{1})+1.$

\item[(b)] $G_{1}^{0}$ has a cycle $C=v_{1},\ldots ,v_{t}v_{1}$ such that $%
v_{1}$ is adjacent to a vertex in $G_{2}^{0}$, say $u$, and $%
N_{G}(v_{t})\subseteq V(G_{2}^{0})\cup V(C)$.\newline
Since $G$ has no cycle of length $\equiv 0\;\mathrm{or}\;2\pmod 3$, we have $%
t\equiv 1\pmod 3$. Let $G_{2}^{1}$ be the graph obtained from $G_{2}^{0}$
and the cycle $C$ by adding the edge $v_{1}u$. By Lemma \ref%
{tailedcycle-3p+1}, $\overrightarrow{f}$ can be extended to a 3-triple of
RDFs $\overrightarrow{g}$ of $G_{2}^{1}$ such that $\omega (g,C)\leq 2n(C)$
and each vertex in $V(C)-\{v_{t}\}$ is $\overrightarrow{g}$-strong. In
addition, it is clear that $\omega (\overrightarrow{g})\leq \omega (%
\overrightarrow{f})+\omega (\overrightarrow{g},C)\leq 2n(G_{2}^{1})+1$.

\item[(c)] $G_{1}^{0}$ contains a tailed $m$-cycle $C_{m,\ell }$, such that $%
y_{\ell }$ is adjacent to some vertex in $G_{2}^{0}$, say $u$, and $%
N_{G}(x_{2})\cup N_{G}(x_{m})\subseteq V(G_{2}^{0})\cup V(C_{m,\ell })$.%
\newline
As above in (b), since $G$ has no cycle of length $\equiv 0\;\mathrm{or}\;2%
\pmod 3$, we have $m\equiv 1\pmod 3$. Let $G_{2}^{1}$ be the graph obtained
from $G_{2}^{0}$ and the tailed $m$-cycle $C_{m,\ell }$ by adding the edge $%
uy_{\ell }$. By Lemma \ref{tailedcycle-3p+1}, $\overrightarrow{f}$ can be
extended to a 3-triple of RDFs $\overrightarrow{g}$ of $G_{2}^{1}$ such that
$\omega (\overrightarrow{g},C_{m,\ell })\leq 2n(C_{m,\ell })$ and each
vertex of $C_{m,\ell }$ but $x_{m}$\ is $\overrightarrow{g}$-strong.
Therefore, we also have $\omega (\overrightarrow{g})\leq 2n(G_{2}^{1})+1.$
\end{description}

Now, let $G_{1}^{1}=G-G_{2}^{1}$. By repeating the above process, we obtain
a $3$-tuple of RDFs $G$ that is $\overrightarrow{h}=(h_{1},h_{2},h_{3})$
such that $\omega (\overrightarrow{h})\leq 2n(G)+1\leq \frac{3n(4k+8)}{6k+8}$%
. Therefore, $\omega (h_{j})\leq \frac{n(4k+8)}{6k+8}$ for some $j\in
\{1,2,3\},$ and this completes the proof.$\hfill \Box $

\begin{theorem}
\label{Th2}\emph{\emph{Let} $k\geq 1$ be an integer and let $G$ be a
connected graph of order $n\geq 6k+9$ with minimum degree at least 2 and
having a cycle $C$ with length $\equiv 0\pmod
3$ such that any other cycle of $G$ with length $\equiv 0\;\mathrm{or}\;2%
\pmod 3$ has at least a common vertex with $C$. Then $\gamma _{R}(G)\leq
\frac{(4k+8)n}{6k+11}$.}
\end{theorem}

\noindent \textbf{Proof.} Assume that the vertices of the cycle $C$ with
length $\equiv 0\pmod
3$ are labelled by $z_{1}z_{2}\ldots z_{r}z_{1}.$ If $V(G)=V(C)$, then
clearly $\gamma _{R}(G)\leq \frac{2n}{3}<\frac{(4k+8)n}{6k+11}$. Hence, we
assume that $V(C)\subsetneqq V(G)$. Let $G_{2}^{0}=C$ and let $G_{1}^{0}$ be
the graph induced by $V(G)-V(G_{2}^{0}).$ By Lemma \ref{cyclepath}, $%
G_{2}^{0}$ has a 3-tuple of RDFs $\overrightarrow{f}=(f_{1},f_{2},f_{3})$
such that $\omega (\overrightarrow{f})\leq 2r$ and all vertices of $C$ are $%
\overrightarrow{f}$-strong. Now, according to Lemma \ref{12}, we consider
the following three possibilities.

\begin{description}
\item[(a)] $G_{1}^{0}$ has a path $P=v_{1}\ldots v_{t}$ such that $%
v_{1},v_{t}$ have neighbors in $V(G_{2})$, say $u,v$ (possibly $u=v$), and $%
N_{G}(v_{1})\cup N_{G}(v_{t})\subseteq V(G_{2}^{0})\cup V(P)$.\newline
Let $G_{2}^{1}$ be the graph obtained from $G_{2}^{0}$ and the path $P$ by
adding the edges $v_{1}u$ and $v_{t}v$. By Lemma \ref{ear1}, $%
\overrightarrow{f}$ can be extended to a 3-triple of RDFs $\overrightarrow{g}
$ of $G_{2}^{1}$ such that $\omega (g,P)\leq 2n(P)=2t$ and each vertex in $%
V(P)-\{v_{1},v_{t}\}$ is $\overrightarrow{g}$-strong. In this case, we have $%
\omega (\overrightarrow{g})=\omega (\overrightarrow{f})+\omega (%
\overrightarrow{g},P)\leq 2n(G_{2}^{1}).$

\item[(b)] $G_{1}^{0}$ has a cycle $C^{\prime }=v_{1},\ldots ,v_{t}v_{1}$
such that $v_{1}$ is adjacent to a vertex in $G_{2}^{0}$, say $u$, and $%
N_{G}(v_{t})\subseteq V(G_{2}^{0})\cup V(C^{\prime })$.\newline
By assumption, we have $t\equiv 1\pmod 3$. Let $G_{2}^{1}$ be\ the graph
obtained from $G_{2}^{0}$ and the cycle $C^{\prime }$ by adding the edge $%
v_{1}u$. By Lemma \ref{tailedcycle-3p+1}, $\overrightarrow{f}$ can be
extended to a 3-triple of RDFs $\overrightarrow{g}$ of $G_{2}^{1}$ such that
$\omega (\overrightarrow{g},C^{\prime })\leq 2n(C^{\prime })=2t$ and each
vertex in $V(C^{\prime })-\{v_{t}\}$ is $\overrightarrow{g}$-strong.
Moreover, we also obtain $\omega (\overrightarrow{g})\leq 2n(G_{2}^{1})$.

\item[(c)] $G_{1}^{0}$ contains a tailed $m$-cycle $C_{m,\ell }$, such that $%
y_{\ell }$ is adjacent to some vertex in $G_{2}^{0}$, say $u$, and $%
N_{G}(x_{2})\cup N_{G}(x_{m})\subseteq V(G_{2}^{0})\cup V(C_{m,\ell })$.%
\newline
As above, $m\equiv 1\pmod 3$. Let $G_{2}^{1}$ be the graph obtained from $%
G_{2}^{0}$ and the tailed $m$-cycle $C_{m,\ell }$ by adding the edge $%
uy_{\ell }$. By Lemmas \ref{tailedcycle-3p+1}, $\overrightarrow{f}$ can be
extended to a 3-triple of RDFs $\overrightarrow{g}$ of $G_{2}^{1}$ such that
such that $\omega (\overrightarrow{g},C_{m,\ell })\leq 2n(C_{m,\ell })$ and
each vertex of $C_{m,\ell }$ but $x_{m}$\ is $\overrightarrow{g}$-strong. In
addition, we have $\omega (\overrightarrow{g})\leq 2n(G_{2}^{1})$.
\end{description}

Let $G_{1}^{1}=G-G_{2}^{1}$. By repeating the above process, we obtain a
3-tuple of RDFs $G$ that is $\overrightarrow{h}=(h_{1},h_{2},h_{3})$ such
that $\omega (\overrightarrow{h})\leq 2n(G)+1\leq \frac{3n(4k+8)}{6k+8}$.
Therefore $\omega (f_{j})\leq \frac{n(4k+8)}{6k+8}$ for some $j\in
\{1,2,3\}, $ and this completes the proof.$\hfill \Box $

\begin{theorem}
\label{Th3}\emph{Let $k\geq 1$ be an integer and let $G$ be a connected
graph of order $n\geq 6k+9\;$and minimum degree at least 2 which does not
contain {neither }any induced $\{C_{5},C_{8},\ldots ,C_{3k+2}\}$-cycles {nor
}any cycle of length $\equiv 0\pmod 3$, and every two distinct cycles of
length $\equiv 2\pmod 3$ (if any) have at least a common vertex. If $G$ has
a cycle $C$ with length $\equiv 2\pmod 3$, then $\gamma _{R}(G)\leq \frac{%
(4k+8)n}{6k+11}$.}
\end{theorem}

\noindent \textbf{Proof.} Let $C=z_{1}z_{2}\ldots z_{p}z_{1}$ be a cycle of
length $\equiv 2\pmod 3$ in $G$ {chosen first not induced, if it exists,
otherwise it is of course induced. }If $V(G)=V(C)$, then we have $\gamma
_{R}(G)\leq \frac{2n+2}{3}<\frac{(4k+8)n}{6k+11}$. Hence, we can assume that
$V(C)\subsetneqq V(G)$.

First assume there is either a cycle $C^{\prime }=x_{1}x_{2}\ldots
x_{m}x_{1} $ such that $x_{1}$ is adjacent to a vertex {of }$C$, say $z_{1}$%
, and $N_{G}(x_{m})\subseteq V(C)\cup V(C^{\prime })$, or a tailed $m$-cycle
$C_{m,\ell }$ in $G$ such that $y_{\ell }$ is adjacent to a vertex {of} $C$,
say $z_{1}$, and $N_{G}(x_{m})\subseteq V(C)\cup V(C_{m,\ell })$. By
assumption $m\equiv 1\pmod 3$. Let $G_{2}^{0}=C+C^{\prime }+x_{1}z_{1}$ or $%
G_{2}^{0}=C+C_{m.\ell }+y_{\ell }z_{1}$ {(depending on which situation
occurs, the first or the second one), }and let $G_{1}^{0}=G-G_{2}^{0}$. By
Lemma \ref{MainLem}, $G_{2}^{0}$ has a 3-tuple of RDFs $\overrightarrow{f}%
=(f_{1},f_{2},f_{3})$ such that $\omega (\overrightarrow{f})\leq
2n(G_{2}^{0})+1$ and all vertices of $G_{2}^{0}$ but $x_{m}$ are $%
\overrightarrow{f}$-strong. Considering our assumption and Lemma \ref{12}, {%
one of the following situations holds. }

\begin{description}
\item[(a)] $G_{1}^{0}$ has a path $P=v_{1},\ldots ,v_{t}$ such that $v_{1}$
and $v_{t}$ are adjacent to some vertices in $V(G_{2}^{0})$, say $u,v$
(possibly $u=v$) and $N_{G}(v_{1})\cup N_{G}(v_{t})\subseteq
V(G_{2}^{0})\cup V(P)$. {We note that }$x_{m}\notin \{u,v\}$, {since }$%
C^{\prime }$ or $C_{m,\ell }$ {has been chosen so that }$x_{m}$ satisfies $%
N_{G}(x_{m})\subseteq V(C)\cup V(C_{m,\ell })$. {Hence }$u$ and $v$ are $%
\overrightarrow{f}$-strong. \newline
Let $G_{2}^{1}$ be {the graph }obtained from $G_{2}^{0}$ and the path $P$ by
adding the edges $v_{1}u,v_{t}v$. By Lemma \ref{ear1}, $\overrightarrow{f}$
can be extended to a 3-triple $\overrightarrow{g}$ such that $\omega (%
\overrightarrow{g})\leq 2n(G_{2}^{1})+1$ and all vertices of $%
V(P)-\{v_{1},v_{t}\}$ are $\overrightarrow{g}$-strong.

\item[(b)] $G_{1}^{0}$ has a cycle $C^{\prime \prime }=v_{1},\ldots
,v_{t}v_{1}$ such that $v_{1}$ is adjacent to a vertex in $G_{2}^{0}$, say $%
u $, and $N_{G}(v_{t})\subseteq V(G_{2}^{0})\cup V(C^{\prime \prime })$. {A
same argument as in item (a) shows that }$u\neq x_{m}$, {and thus }$u$ is $%
\overrightarrow{f}$-strong. \newline
By assumption, we have $t\equiv 1\pmod 3$. Let $G_{2}^{1}$ be {\ the graph }%
obtained from $G_{2}^{0}$ and the cycle $C^{\prime \prime }$ by adding the
edge $v_{1}u$. By Lemma \ref{tailedcycle-3p+1}, $\overrightarrow{f}$ can be
extended to a 3-triple $\overrightarrow{g}$ such that $\omega (%
\overrightarrow{g})\leq 2n(G_{2}^{1})+1$ and each vertex $V(C^{\prime \prime
})-\{v_{t}\}$ is $\overrightarrow{g}$-strong.

\item[(c)] $G_{1}^{0}$ contains a tailed $m^{\prime }$-cycle $C_{m^{\prime
},\ell ^{\prime }}$, such that $y_{\ell ^{\prime }}$ is adjacent to some
vertex in $G_{2}^{0}$, say $u$, and $N_{G}(x_{2})\cup N_{G}(x_{m})\subseteq
V(G_{2}^{0})\cup V(C_{m^{\prime },\ell ^{\prime }})$. {Note that }$u\neq
x_{m}$ and $u$ is $\overrightarrow{f}$-strong. \newline
As above $m^{\prime }\equiv 1\pmod 3$. Let $G_{2}^{1}$ be {the graph }%
obtained from $G_{2}^{0}$ and the tailed cycle $C_{m^{\prime },\ell ^{\prime
}}$ by adding the edge $uy_{\ell ^{\prime }}$. By Lemma \ref%
{tailedcycle-3p+1}, $\overrightarrow{f}$ can be extended to a 3-triple $%
\overrightarrow{g}$ such that $\omega (\overrightarrow{g})\leq
2n(G_{2}^{1})+1$ and all vertices of $V(C_{m^{\prime },\ell ^{\prime
}})-\{x_{m^{\prime }}\}$ are $\overrightarrow{g}$-strong.
\end{description}

Let $G_{1}^{1}=G-G_{2}^{1}$. By repeating {the }above process, we obtain a
3-tuple of RDFs $\overrightarrow{g^{\prime }}=(g_{1}^{\prime },g_{2}^{\prime
},g_{3}^{\prime })$ such that $\omega (\overrightarrow{g^{\prime }})\leq
2n(G)+1\leq \frac{3n(4k+8)}{6k+8}$. It follows that $\omega (g_{j}^{\prime
})\leq \frac{n(4k+8)}{6k+8}$ for some $j\in \{1,2,3\}$ as desired.

{Next we can }assume that there is {neither a }cycle $C^{\prime
}=(x_{1}x_{2}\ldots x_{m}x_{1})$ such that $x_{1}$ is adjacent to a vertex
in $C$ and $N_{G}(x_{m})\subseteq V(C)\cup V(C^{\prime })$, {nor a }tailed $%
m $-cycle $C_{m,\ell }$ with $m\equiv 1\pmod 3$ in $G$ such that $y_{\ell }$
is adjacent to a vertex in $C$ and $N_{G}(x_{m})\subseteq V(C)\cup
V(C_{m,\ell })$. Let $H_{2}^{0}=C$ and $H_{1}^{0}=G-H_{2}^{0}$. It follows
from Lemma \ref{12} and the assumptions that $H_{1}^{0}$ has a path $%
P=v_{1},\ldots ,v_{t}$ such that $v_{1},v_{t}$ are adjacent to some vertices
in $V(H_{2}^{0})$, say $z_{1},z_{j}$ (possibly $j=1$) and $N_{G}(v_{1})\cup
N_{G}(v_{t})\subseteq V(H_{2}^{0})\cup V(P)$. We consider the following
cases.

\smallskip \noindent \textbf{Case 1.} $j=1$.\newline
Let $G_{2}^{1}$ be {the graph }obtained from $H_{2}^{0}$ and the path $P$ by
adding the edges $v_{1}z_{1},v_{t}z_{1}$ and let $G_{1}^{1}=G-G_{2}^{1}$. By
Lemma \ref{MainLem}{-(1)}, $G_{2}^{1}$ has a triple $\overrightarrow{g}$ of
RDFs such that $\omega (\overrightarrow{f})\leq 2n(G_{2}^{1})+1$ and all
vertices of $G_{2}^{1}$ but $v_{1},v_{t}$ are $\overrightarrow{g}$-strong.
If $V(G)=V(G_{2}^{1})$ {(and hence }$G_{1}^{1}$ is empty), then the result
follows. {Hence, a}ssume that $V(G)\neq V(G_{2}^{1})$. By the assumptions
and Lemma \ref{12}, we deduce that $G_{1}^{1}$ has a path $P^{\prime
}=v_{1}^{\prime },\ldots ,v_{t^{\prime }}^{\prime }$ such that $%
v_{1}^{\prime },v_{t^{\prime }}^{\prime }$ are adjacent to some vertices in $%
V(G_{2}^{1})$, say $u,v$ (possibly $u=v$) and $N_{G}(v_{1}^{\prime })\cup
N_{G}(v_{t^{\prime }}^{\prime })\subseteq V(G_{2}^{1})\cup V(P^{\prime })$.
Let $G_{2}^{2}$ be obtained from $G_{2}^{1}$ and the path $P^{\prime }$ by
adding the edges $v_{1}^{\prime }u,v_{t^{\prime }}^{\prime }v$ and let $%
G_{1}^{2}=G-G_{2}^{2}$. {Note that }$v_{1},v_{t}\notin \{u,v\}$ {and thus }$%
u,v$ are $\overrightarrow{g}$-strong. By Lemma \ref{ear1}, $\overrightarrow{g%
}$ can be extended to a 3-triple $\overrightarrow{g^{\prime }}$ such that $%
\omega (\overrightarrow{f_{2}})\leq 2n(G_{2}^{2})+1$ and all vertices of $%
V(P^{\prime })-\{v_{1}^{\prime },v_{t^{\prime }}^{\prime }\}$ are $%
\overrightarrow{g^{\prime }}$-strong. By repeating {the }above process, we
obtain a 3-tuple of RDFs $\overrightarrow{g^{\ast }}=(g_{1}^{\ast
},g_{2}^{\ast },g_{3}^{\ast })$ such that $\omega (\overrightarrow{g^{\ast }}%
)\leq 2n(G)+1\leq \frac{3n(4k+8)}{6k+8}$. It follows that $\omega
(g_{r}^{\ast })\leq \frac{n(4k+8)}{6k+8}$ for some $r\in \{1,2,3\}$ as
desired.

\smallskip \noindent \textbf{Case 2.} $j\neq 1$.\newline
We distinguish the following {three }subcases.

\smallskip \textbf{Subcase 2.1.} $t\equiv 1\pmod 3$.\newline
Since $G$ has no cycle of length $\equiv 0\pmod 3$, we have $j\not\equiv 2%
\pmod 3$. Let $G_{2}^{1}$ be {the graph }obtained from $H_{2}^{0} $ and the
path $P$ by adding the edges $v_{1}z_{1},v_{t}z_{j}$ and let $%
G_{1}^{1}=G-G_{2}^{1}$. By Lemma \ref{ear2}, $G_{2}^{1}$ has a triple $%
\overrightarrow{f}$ of RDFs such that $\omega (\overrightarrow{f})\leq
2n(G_{2}^{1})+1$ and all vertices of $G_{2}^{1}$ but $v_{1},v_{t}$ are $%
\overrightarrow{f_{1}}$-strong. As in Case 1, we can obtain a 3-tuple of
RDFs $\overrightarrow{g^{\ast }}=(g_{1}^{\ast },g_{2}^{\ast },g_{3}^{\ast })$
such that $\omega (\overrightarrow{g^{\ast }})\leq 2n(G)+1\leq \frac{3n(4k+8)%
}{6k+8}$ {yielding the desired result}.

\smallskip \textbf{Subcase 2.2.} $t\equiv 2\pmod 3$.\newline
{Observe that if }$j\equiv 0\pmod 3$, {then }$%
z_{1}v_{1}...v_{t}z_{j}z_{j+1}...z_{p}z_{1}$ would be a cycle of length $%
\equiv 0\pmod 3$, a contradiction, {and if $j\equiv 1\pmod 3$, {%
then }$z_{1}v_{1}...v_{t}z_{j}z_{j-1}...z_{2}z_{1}$ would be a cycle of
length $\equiv 0\pmod 3$, a contradiction again. Hence $j\equiv 2\pmod 3$}. {%
Now, as }in Subcase 2.1, we can get the result.

Considering Subcases 2.1 and 2.2, we may assume all ears of $C$ in $%
G_{1}^{0} $ have length $\equiv 0\pmod 3$.

\smallskip \textbf{Subcase 2.3.} $t\equiv 0\pmod 3$.\newline
Considering the cycles generated by $C+P+v_{1}z_{1}+v_{t}z_{j}$ and that
fact that $G$ has no cycle of length $\equiv 0\pmod 3$, we {\ deduce that }$%
j\equiv 2\pmod 3$. Let $C_{1}=(z_{1}z_{2}\ldots z_{p}z_{1})$, $%
C_{2}=(z_{1}v_{1}v_{2}\ldots v_{t}z_{j}z_{j-1}\ldots z_{2}z_{1})$ and $%
C_{3}=(z_{1}v_{1}v_{2}\ldots v_{t}z_{j}z_{j+1}\ldots z_{p}z_{1})$. Clearly
the cycles $C_{1},C_{2}$, $C_{3}$ are {all }of length $\equiv 2\pmod 3$.

{Assume first that }$C_{1}$ is not an induced cycle in $G$. Then by Lemma %
\ref{induced} and considering the ear we have $n(C_{1}\cup C_{2})\geq 6k+11$%
. Let $G_{2}^{0}=C_{1}\cup C_{2}$ and $G_{1}^{0}=G-G_{2}^{0} $. It is not
hard to see that $G_{2}^{0}$ has a 3-tuple $\overrightarrow{f}$ of RDFs such
that $\omega (\overrightarrow{f})\leq 2n(G_{2}^{0})+2$ and all vertices of $%
G_{2}^{0}$ but $v_{1},v_{t}$ are $\overrightarrow{f}$-strong. If $%
V(G)=V(G_{2}^{0})$, then the result follows. {Hence }assume that $V(G)\neq
V(G_{2}^{0})$. By the assumptions and Lemma \ref{12}, we deduce that $%
G_{1}^{0}$ has a path $P^{\prime }=v_{1}^{\prime },\ldots ,v_{t^{\prime
}}^{\prime }$ such that $v_{1}^{\prime },v_{t^{\prime }}^{\prime }$ are
adjacent to some vertices in $V(G_{2}^{0})$, say $u,v$ (possibly $u=v$) and $%
N_{G}(v_{1}^{\prime })\cup N_{G}(v_{t^{\prime }}^{\prime })\subseteq
V(G_{2}^{0})\cup V(P^{\prime })$. {Note that }$u,v\notin \{v_{1},v_{t}\}$ {%
since }$N_{G}(v_{1})\cup N_{G}(v_{t})\subseteq V(H_{2}^{0})\cup V(P)$. Thus $%
{u,v}$ are $\overrightarrow{f}$-strong. {Now, let }$G_{2}^{1}$ be {the graph
}obtained from $G_{2}^{1}$ and the path $P^{\prime }$ by adding the edges $%
v_{1}^{\prime }u,v_{t^{\prime }}^{\prime }v$ and let $G_{1}^{1}=G-G_{2}^{1}$%
. By Lemma \ref{ear1}, $\overrightarrow{f} $ can be extended to a 3-triple $%
\overrightarrow{g}$ such that $\omega (\overrightarrow{g})\leq
2n(G_{2}^{2})+1$ and all vertices\ of $P^{\prime }$ but $v_{1}^{\prime
},v_{t^{\prime }}^{\prime }$ are $\overrightarrow{f_{2}}$-strong. By
repeating above process, we obtain a 3-tuple of RDFs $\overrightarrow{%
g^{\ast }}=(g_{1}^{\ast },g_{2}^{\ast },g_{3}^{\ast })$ such that $\omega (%
\overrightarrow{g^{\ast }})\leq {2n(G)+2}\leq \frac{3n(4k+8)}{%
6k+8}$. It follows that $\omega (\overrightarrow{g_{j}^{\ast }})\leq \frac{%
n(4k+8)}{6k+8}$ for some $j\in \{1,2,3\}$ as desired.

Assume now that $C_{1}$ is {an induced cycle}. By the choice of $C,$ we may
assume that $G$ has no cycle of length $\equiv 2\pmod 3$ which is not
induced. {Hence the cycle }$C_{2}$ {is also induced. }Let $%
G_{2}^{0}=C_{1}\cup C_{2}$ and $G_{2}^{1}=G-G_{2}^{0}$. {There are the
following two possibilities. }

\begin{itemize}
\item $V(G)=V(G_{2}^{0})$.\newline
Suppose $n(C_{1})=3t_{1}+2$ and $t=3t_{2}$. Using the fact that $n\geq 6k+9,$
we obtain
\begin{equation*}
\begin{array}{lll}
n & = & n(C_{1})+t \\
& = & 3t_{1}+2+3t_{2} \\
& \geq & 3(2k+3)%
\end{array}%
\end{equation*}%
implying that $t_{1}+t_{2}\geq 2k+3-2/3$. {Since }$t_{1}+t_{2}$ {is integer,
we deduce that }$t_{1}+t_{2}\geq 2k+3,$ and thus $n\geq 6k+11$. {Now, it }is
easy to see that $\gamma _{R}(G)\leq \frac{2n+2}{3}\leq \frac{(4k+8)n}{6k+11}
$.

\item $V(G_{2}^{0})\subsetneqq V(G)$.\newline
Clearly $G_{2}^{0}$ has a triple $\overrightarrow{f^{0}}$ of RDFs such that $%
\omega (\overrightarrow{f^{0}})\leq n(G_{2}^{0})+2$ and all vertices of $%
G_{2}^{0}$ but $v_{1},v_{t}$ are $\overrightarrow{f^{0}}$-strong. By the
assumptions and Lemma \ref{12}, we deduce that $G_{1}^{0}$ has a path $%
P_{1}=v_{1}^{1},\ldots ,v_{q_{1}}^{1}$ such that $v_{1}^{1},v_{q_{1}}^{1}$
are adjacent to some vertices in $V(G_{2}^{0})$, say $u,v$ (possibly $u=v$)
and $N_{G}(v_{1}^{1})\cup N_{G}(v_{q_{1}}^{1})\subseteq V(G_{2}^{1})\cup
V(P_{1})$. {Recall that }$u,v\in \{v_{1},v_{t}\}$ {\ and thus they are }$%
\overrightarrow{f^{0}}$-strong. {Moreover, }since every cycle of $G$
intersects $C_{1}$, we have $V(C_{1})\cap \{u,v\}\neq \emptyset $. {Hence
vertices }$u,v$ {may belong to }$C_{1}$, $C_{2}$ or $C_{3}$. {Now, }seeing
Case 1 and Subcase 2.1 and 2.2, we may assume that $q_{1}\equiv 0\pmod 3$.
Let $q_{1}=3q_{1}^{\prime }$ and let $G_{2}^{1}$ be {the graph }obtained
from $G_{2}^{0}$ and the path $P_{1}$ by adding the edges $%
v_{1}^{1}u,v_{q_{1}}^{1}v$ and let $G_{1}^{1}=G-G_{2}^{1}$. By Lemma \ref%
{ear1}, $\overrightarrow{f^{0}}$ can be extended to a 3-triple $%
\overrightarrow{f^{1}}$ {of }$G_{2}^{1}$ such that {$\omega (\overrightarrow{%
f^{1}})\leq 2n(G_{2}^{1})+2$} and all vertices {\ of }$P_{1}$ but $%
v_{1}^{1},v_{q_{1}}^{1}$ are $\overrightarrow{f^{1}}$-strong. If $%
V(G)=V(G_{2}^{1})$, then $n=3t_{1}+3t_{2}+3q_{1}^{\prime }+2$. As above we
can see that $n\geq 6k+11,$ {implying that }$\gamma _{R}(G)\leq \frac{2n+2}{3%
}\leq \frac{(4k+8)n}{6k+11}$. {Hence assume that }$V(G_{2}^{1})\subsetneqq
V(G)$. By the assumptions and Lemma \ref{12}, we deduce that $G_{1}^{1}$ has
a path $P_{2}=v_{1}^{2},\ldots ,v_{q_{2}}^{2}$ such that $%
v_{1}^{2},v_{q_{2}}^{2}$ are adjacent to some vertices in $V(G_{2}^{1})$,
say $u^{\prime },v^{\prime }$ (possibly $u^{\prime }=v^{\prime }$) and $%
N_{G}(v_{1}^{2})\cup N_{G}(v_{q_{2}}^{2})\subseteq V(G_{2}^{1})\cup V(P_{2})$%
. Since every cycle of $G$ intersects $C_{1}$, we have $V(C_{1})\cap
\{u^{\prime },v^{\prime }\}\neq \emptyset $. On the other hand, we note that
$u^{\prime },v^{\prime }$ lies on a cycle of length $\equiv 2\pmod 3$.
Seeing Case 1 and Subcase 2.1 and 2.2, we may assume that $q_{2}\equiv 0%
\pmod 3$. Let $q_{2}=3q_{2}^{\prime }$ and let $G_{2}^{2}$ be {the graph }%
obtained from $G_{2}^{1}$ and the path $P_{2}$ by adding the edges $%
v_{1}^{2}u,v_{q_{2}}^{2}v$ and let $G_{1}^{2}=G-G_{2}^{2}$. By Lemma \ref%
{ear1}, $\overrightarrow{f^{1}}$ can be extended to a 3-triple $%
\overrightarrow{f^{2}}$ such that $\omega (\overrightarrow{f^{2}})\leq
2n(G_{2}^{2})+2$ and all {of }$P_{2}$ but $v_{1}^{2},v_{q_{2}}^{2}$ are $%
\overrightarrow{f^{2}}$-strong. If $V(G)=V(G_{2}^{2})$, then $%
n=3t_{1}+3t_{2}+3q_{1}^{\prime }+3q_{2}^{\prime }+2$. As above we can see
that $n\geq 6k+11,$ {implying that }$\gamma _{R}(G)\leq \frac{2n+2}{3}\leq
\frac{(4k+8)n}{6k+11}$. {\ Hence suppose that }$V(G_{2}^{k})\subsetneqq V(G)$%
. ~ By repeating {the }above process, we obtain a subgraph $G_{2}^{k}$ with $%
n(G_{2}^{k})\geq 6k+11$ and having a 3-tuple $\overrightarrow{f^{k}}$ of
RDFs such that $\omega (\overrightarrow{f^{k}})\leq {
2n(G_{k}^{2})+2}$ and all vertices of $G_{2}^{k}$ with a neighbor outside of
$G_{2}^{k}$ are $\overrightarrow{f^{k}}$-strong. If $V(G)=V(G_{2}^{k})$,
then the result follows immediately. {Otherwise, let $%
G_{1}^{k}=G-V(G_{2}^{k})$}. Now using Lemma \ref{ear1} we can obtain extend $%
G_{2}^{k}$ to a subgraph $G_{2}^{k+1}$ by adding an ear in $G_{1}^{k}$ and
extend $\overrightarrow{f^{k}}$ to a 3-tuple $\overrightarrow{f^{k+1}}$ of
RDFs such that $\omega (\overrightarrow{f^{k+1}})\leq 2n(G_{k+1}^{2})+2$ and
all vertices of $G_{2}^{k+1}$ with a neighbor outside of $G_{2}^{k+1}$ are $%
\overrightarrow{f^{k+1}}$-strong. By repeating this process we obtain a
3-tuple $g$ of RDFs of $G$ such that $\omega (\overrightarrow{g})\leq
2n(G)+2 $ and this leads to the result as above. $\hfill \Box $
\end{itemize}

\section{Some more lemmas}

Let $\mathcal{F}_{i}$ be the family of all cycles of length $\equiv i\pmod 3$
with $i\in \{0,1,2\}$. Let $\mathcal{F}_{0,2}$ be the family of all {%
connected }graphs obtained from a cycle $C$ of $\mathcal{F}_{0}$ and a cycle
$C^{\prime }$ of $\mathcal{F}_{2}$ by {joining a vertex }$x$ {of }$C$ {\and %
a vertex }$y$ {of }$C^{\prime }$ {by either an edge }$xy${\ or by a
nontrivial path that we add so that one of the envertices of the path is
attached to }$x$ and the other one to $y$; $\mathcal{F}_{2,2}$ be the family
of all {connected} graphs obtained from two cycles in $\mathcal{F}_{2}$ by
adding an edge between them; and let $\mathcal{F}_{3}$ be the family of all
graphs $G$ obtained from a graph $G^{\prime }$ in $\mathcal{F}_{0,2}$ and a
graph $G^{\prime \prime }$ in $\mathcal{F}_{2,2}$ by adding {either an edge
or a path joining a vertex of }$G^{\prime }$ to a vertex of $G^{\prime
\prime }$ {so that all vertices of the path become of degree two in }$G.$

\bigskip\

Let $\mathcal{B}_{r,s}\;(r+s\geq 2)$ be the family of {connected }graphs
obtained from $r$ tailed cycles $C_{n_{1},\ell _{1}},\ldots ,C_{n_{r},\ell
_{r}}$ and $s$ cycles $C_{m_{1}},\ldots ,C_{m_{s}}$, where $n_{i}\equiv 2%
\pmod 3$ and $m_{j}\equiv 2\pmod 3$ for each $i,j$, by adding a new vertex $%
z $ ({which we call }special vertex) and {\ joining by edges }$z$ to the {%
unique leaf }of each graph $C_{n_{i},\ell _{i}}$ and to one vertex of each
cycle $C_{m_{j}}$. {Moreover, each of the }$s$ {cycles will be called a
\textit{near cycle} of }$z.${\ }Let $\mathcal{E}=\cup _{r,s\geq 0;r+s\geq 2}%
\mathcal{B}_{r,s}$.

\begin{lemma}
\label{Lam3}\emph{Let $G$ be a connected graph {with at least two disjoint
cycles and let }$\mathcal{F}$ be a family of pairwise disjoint cycles of
length $\equiv 0,2\pmod 3$ in $G$ with $|\mathcal{F}|\geq 2$. Then $G$ has
two disjoint subgraphs $G_{1}$ (possibly null) and $G_{2}$ such that $%
V(G)=V(G_{1})\cup V(G_{2})$, $G_{1}$ has no cycle of $\mathcal{F}$ and each
component of $G_{2}$ is in $\mathcal{F}_{0}\cup \mathcal{F}_{0,2}\cup
\mathcal{F}_{2,2}\cup \mathcal{E}.$}
\end{lemma}

\noindent \textbf{Proof.} The proof is by induction on the number of cycles
in $\mathcal{F}$. First let $|\mathcal{F}|=2$ with $\mathcal{F}%
=\{C_{1},C_{2}\}$. {Since }$G$ {is connected, let }$P $ be a shortest path {%
joining a vertex of }$C_{1}$ {to a vertex of $C_{2}$. If both $C_{1},C_{2}$
have length $\equiv 0\pmod 3$, then let $G_{2}=C_{1}\cup C_{2}$ and if one
of the {two }cycles has length $\equiv 2\pmod 3$, then let $%
G_{2}=C_{1}+C_{2}+P$. Assume that $G_{1}=G-V(G_{1})$. Clearly $G_{1}$ has no
cycle of $\mathcal{F}$ and each component of $G_{2}$ belongs to $\mathcal{F}%
_{0}\cup \mathcal{F}_{0,2}\cup \mathcal{F}_{2,2}\cup \mathcal{E}$, {%
establishing the base case.} }

\noindent Next let $|\mathcal{F}|=3$ and $\mathcal{F}=\{C_{1},C_{2},C_{3}\}$%
. Assume that $C_{i}=x_{1}^{i}x_{2}^{i}\ldots x_{n_{i}}^{i}x_{1}^{i}$ for $%
i\in \{1,2,3\}$. If {each cycle of }$\mathcal{F}$ has length $\equiv 0%
\pmod 3
$, then let $G_{2}=C_{1}\cup C_{2}\cup C_{3}$ and $G_{1}=G-G_{2}$. {Clearly
the result holds. Hence assume that }one of the {three }cycles has length $%
\equiv 2\pmod 3$, say $C_{1}$. Let $P$ be a shortest path {joining a vertex
of }$C_{1}$ {to a vertex in }$C_{2}$ or $C_{3}$. Assume, without loss of
generality, that $P$ {joins }$C_{1}$ and $C_{2}$, where $%
P=(x_{1}^{1}=)z_{0}z_{1}...z_{k}(=x_{1}^{2}).$ If $C_{3}$ has length $\equiv
0\pmod 3,$ {then by setting }$G_{2}=\left( C_{1}\cup C_{2}\cup P\right) \cup
C_{3}$ and $G_{1}=G-G_{2},$ {it is clear that the result holds. Hence we
assume that $C_{3}$ has length $\equiv 2\pmod 3.$ {Now, let }$%
Q=(x_{1}^{3}=)y_{0}y_{1}...y_{s}$ be a shortest path {joining a vertex of }$%
C_{3}$ {to a vertex }$y_{s}$ {belonging to }$V(C_{1})\cup V(C_{2})\cup V(P)$%
. Assume that $y_{s}\in V(P)-\{x_{1}^{1},x_{1}^{2}\}$, say $y_{s}=z_{m}.$ If
$C_{2}$ has length $\equiv 0\pmod 3,$ then {by setting }$G_{2}=\left(
C_{1}\cup C_{3}\cup P^{\prime }\right) \cup C_{2},$ with $P^{\prime
}=(x_{1}^{1}=)z_{0}z_{1}...z_{m},y_{s-1},\ldots,y_{0},$ and $G_{1}=G-G_{2},$
we get the desired result. Hence we assume that $C_{2}$ has length $\equiv 2%
\pmod 3.$ In this case, the result holds by letting $G_{2}=C_{1}\cup
C_{2}\cup C_{3}\cup P\cup Q$ and $G_{1}=G-G_{2}.$ Finally, assume, without
loss of generality, that\ $y_{s}\in C_{2},$ say $y_{s}=x_{j}^{2}$ (possibly $%
j=1$). Let $G_{2}=C_{1}Px_{1}^{2}\ldots x_{j}^{2}QC_{3}$ and $G_{1}=G-G_{2}$%
. {Note that }$G_{2}$belongs to {$\mathcal{B}_{0,2}\cup \mathcal{B}_{1,1},$}
{and clearly the desired result holds. } }

\noindent Assume now that $|\mathcal{F}|\geq 4$. If all cycles in $\mathcal{F%
}$ have length $\equiv 0\pmod 3$, then the subgraphs $G_{2}=\cup _{i=1}^{|%
\mathcal{F}|}C_{i}$ and $G_{1}=G-G_{2}$ satisfy the conditions {\ and the
result holds. }Hence we assume that one of the cycles in $\mathcal{F} $, {%
say }$C_{0},${\ }has length $\equiv 2\pmod 3$. Let $\mathcal{F}^{\prime }=%
\mathcal{F}-\{C_{0}\}$ and let $G^{1}=G-V(C_{0})$. We consider two cases.

\smallskip \noindent \textbf{Case 1.} $G^{1}$ is connected.\newline
Then $\mathcal{F}^{\prime }$ is a family of disjoint cycles {of length }$%
\equiv 0,2\pmod 3$ in $G^{1}$ {with }\emph{$|\mathcal{F}^{\prime }|\geq 2$}.
By the induction hypothesis, $G^{1}$ has two disjoint subgraphs $%
G_{1}^{\prime }$ (possibly null) and $G_{2}^{\prime }$ such that $%
V(G^{1})=V(G_{1}^{\prime })\cup V(G_{2}^{\prime })$, $G_{1}^{\prime }$ has
no cycle of $\mathcal{F}^{\prime }$, where each component of $G_{2}^{\prime
} $ is in $\mathcal{F}_{0}\cup \mathcal{F}_{0,2}\cup \mathcal{F}_{2,2}\cup
\mathcal{E}.$ Let $H_{1},\ldots ,H_{p}$ be the components of $G_{2}^{\prime
} $. {Suppose without loss of generality that $P:=(x_{1}^{0}=)v_{0}v_{1}%
\ldots v_{t}$ is a shortest path between $V(C_{0})$ and $V(G_{2}^{\prime })$
in $G$ where $v_{t}\in V(G_{2}^{\prime })$.} Without loss of generality,
assume that $v_{t}\in V(H_{1}).$ If $H_{1}\in \mathcal{F}_{0}$, then let $%
H_{1}^{\prime }=H_{1}+P+C_{0}$ and clearly {the }two subgraphs $%
G_{2}=H_{1}^{\prime }\cup H_{2}\cup \ldots H_{p}$ and $G_{1}=G-V(G_{2})$
satisfy the conditions {and result follows.
For the next, we can assume that $H_{1}$ contains at least two cycles. We
distinguish the following. }

\smallskip \textbf{Subcase 1.1.} $H_{1}\in \mathcal{F}_{2,2}$.\newline
Suppose $H_{1}$ is obtained from two cycles $C_{1}=x_{1}^{1}\ldots
x_{m_{1}}^{1}x_{1}^{1}$ and $C_{2}=x_{1}^{2}\ldots x_{m_{2}}^{2}x_{1}^{2}$
by adding a path $Q=(x_{1}^{1}=)z_{0}z_{1}\ldots z_{t}(=x_{1}^{2})$. {We
further assume, without loss of generality, that }$v_{t}=x_{j}^{2}\in
V(C_{2})$. {Let }$H_{1}^{\prime }$ {be the graph obtained from }$C_{0}\cup
C_{1}$ {to which we add the path }$Px_{j-1}^{2}\ldots x_{2}^{2}Q$, {in other
words, }$H_{1}^{\prime }$ {is obtained from }$C_{0}+P+H_{1} $ {by removing
vertices }$x_{j+1}^{2},...,x_{m_{2}}^{2}.$ {Note that }$H_{1}^{\prime }$ {%
belongs to either $\mathcal{F}_{0,2}$ or {$\mathcal{B}_{0,2}\cup\mathcal{B}%
_{1,1} .$ } Now let $G_{2}^{\prime \prime }=H_{1}^{\prime }\cup H_{2}\cup
\ldots H_{p}$ and $G_{1}^{\prime \prime }=G-G_{2}^{\prime \prime }. $ Then
the subgraphs $G_{2}^{\prime \prime }$ and $G_{1}^{\prime \prime }$ satisfy
the conditions {and result follows}. }

\smallskip \textbf{Subcase 1.2.} $H_{1}\in \mathcal{F}_{0,2}$.%
\newline
Using an argument similar to that described in the case $|\mathcal{F}|=3$,
we can obtain two subgraphs $G_{2}^{\prime \prime }$ and $G_{1}^{\prime
\prime }$ satisfying the conditions {and yielding the desired result}.

\smallskip \textbf{Subcase 1.3.} $H_{1}\in \mathcal{B}_{r,s}$ where $r+s\geq
2$.\newline
Let $z^{\ast }$ {be the special vertex of }$H_{1}.$ {If }$v_{t}=z^{\ast }$,
then $H_{1}^{\prime }=H_{1}+P+C_{0}$ is a subgraph belonging to $\mathcal{E}$
and {thus }the subgraphs $G_{2}=H_{1}^{\prime }\cup H_{2}\cup \ldots H_{p}$
and $G_{1}=G-V(G_{2})$ satisfy the conditions and the result follows. Hence
we assume that $v_{t}\neq z^{\ast }.$ First let $r+s=2.$ Then $H_{1}$ is
obtained from two cycles $C_{1}=x_{1}^{1}\ldots x_{m_{1}}^{1}x_{1}^{1}$ and $%
C_{2}=x_{1}^{2}\ldots x_{m_{2}}^{2}x_{1}^{2}$ by adding a path $%
Q=(x_{1}^{1}=)z_{0}z_{1}\ldots z_{t}(=x_{1}^{2}),$where $t\geq 2.$ If $%
v_{t}\in \{z_{1},\ldots ,z_{t-1}\}$, then let $H_{1}^{\prime
}=H_{1}+C_{0}+P. $ {Clearly, }$H_{1}^{\prime }\in \mathcal{B}_{r,s}$ where $%
r+s=3$, and thus the subgraphs $G_{2}=H_{1}^{\prime }\cup H_{2}\cup \ldots
H_{p}$ and $G_{1}=G-V(G_{2})$ satisfy the conditions {and the result yields.
Now, suppose, without loss of generality, that }$v_{t}=x_{j}^{2}\in
V(C_{2}). $ {Let }$H_{1}$ be obtained $C_{0}\cup C_{1}$ {to which we add the
path }$Px_{j-1}^{2}\ldots x_{2}^{2}Q,$ {and set }$G_{2}^{\prime \prime
}=H_{1}^{\prime }\cup H_{2}\cup \ldots H_{p}$ and $G_{1}^{\prime \prime
}=G-G_{2}^{\prime \prime }.$ {Clearly, }$G_{1}^{\prime \prime }$ and $%
G_{2}^{\prime \prime }$ {satisfy the condition and the desired result
follows.}

{Now let } $r+s\geq 3.$ {Assume that }$v_{t}$ {belongs to one of the }${s+r}$
{cycles of }$H_{1},$ {say} $C^{\prime }.$ {\ Let }$H_{1}^{\prime }=C^{\prime
}+C_{0}+P$ and $H_{1}^{\prime \prime }$ be the graph obtained from $H_{1}$
by deleting {the }vertices of $V(C^{\prime })$ and {the path \ (if any)
joining }$z^{\ast }$ to $V(C^{\prime })$ in $H_{1}$. {\ Note that }$%
H_{1}^{\prime }$ belongs to either {$\mathcal{F}_{2,2}$} or {$%
\mathcal{B}_{0,2}\cup \mathcal{B}_{1,1}.$} {Now the subgraphs }$%
G_{2}=H_{1}^{\prime }\cup H_{1}^{\prime \prime }\cup H_{2}\cup \ldots H_{p}$
and $G_{1}=G-V(G_{2})$ satisfy the conditions and {the desired result holds.
Assume that }$v_{t}$ belongs to a path on a tailed cycle $C_{m,\ell }$ of $%
H_{1},$ and let $P^{\prime }$ be the subpath between $v_{t}$ and the cycle $%
C^{\prime }$ of $C_{m,\ell }.$ {Let }$H_{1}^{\prime }=C^{\prime
}+C_{0}+P+P^{\prime }$ and $H_{1}^{\prime \prime }$ be the graph obtained
from $H_{1}$ by deleting {the }vertices of $V(C_{m,\ell }).$ {Note that }$%
H_{1}^{\prime }$ belongs to either {$\mathcal{F}_{2,2}$} or $%
\mathcal{B}_{0,2}.$ {Now the subgraphs }$G_{2}=H_{1}^{\prime }\cup
H_{1}^{\prime \prime }\cup H_{2}\cup \ldots H_{p}$ and $G_{1}=G-V(G_{2})$
satisfy the conditions and {the desired result holds.}

\smallskip \noindent \textbf{Case 2.} $G^{1}$ is disconnected.\newline
Let $M_{1},\ldots ,M_{t}$ be the components of $G^{1}$. Assume first that a
component $M_{i}$ contains all cycles of $\mathcal{F}-\{C_{0}\},$ say $M_{1}$%
. Let $K$ be the subgraph of $G$ induced by $V(C_{0})\cup V(M_{1})$. Clearly
$K$ is connected. Using an argument similar to that described in Case 1 on $%
K-C_{0}$, we get the result.
Henceforth, we may assume that no $M_{i}$ contains all cycles of $\mathcal{F}%
-\{C_{0}\}$ for each $i$. Now, assume that a component $M_{i}$ contains at
least two cycles of $\mathcal{F}$, say $M_{1}.$ Let $G^{2}=G-V(M_{1})$.
Clearly $G^{2}$ is connected. Let $\mathcal{F}_{1}=\{C\mid C\in \mathcal{F}\;%
\mathrm{and}\;V(C)\subseteq V(M_{1})\}$ and $\mathcal{F}_{2}=\mathcal{F}-%
\mathcal{F}_{1}$. By the induction hypothesis, $M_{1}$ has two subgraphs $%
K^{1},K^{2}$ such that $K^{1}$ does not contain any cycle of $\mathcal{F}%
_{1} $ and each component of $K^{2}$ belongs to $\mathcal{F}_{0}\cup
\mathcal{F}_{0,2}\cup \mathcal{F}_{2,2}\cup \mathcal{E}.$ Moreover, $G^{2}$
has two subgraphs $K_{1}^{\prime }$ and $K_{2}^{\prime }$ such that $%
K_{1}^{\prime }$ does not contain any cycle of $\mathcal{F}_{2}$ and each
component of $K_{2}^{\prime } $ belongs to $\mathcal{F}_{0}\cup \mathcal{F}%
_{0,2}\cup \mathcal{F}_{2,2}\cup \mathcal{E}.$\textbf{\ }Now the two
subgraphs $G_{1}=K^{1}\cup K_{1}^{\prime }$ and $G_{2}=K^{2}\cup
K_{2}^{\prime }$ satisfies the conditions yielding the desired result.

{From now on, we can assume that }each $M_{i}$ contains at most one cycle of
$\mathcal{F}$. {Suppose that only the }$s$ first $M_{i}$ contains exactly
one cycle $C_{i}$ of $\mathcal{F}.$ Let {$C_{i}=x_{1}^{i}x_{2}^{i}\ldots
x_{n_{i}}^{i}x_{1}^{i}$ for $0\leq i\leq s.$ In addition, let $%
P_{i}:=(x_{1}^{i}=)w_{0}^{i}\ldots w_{\ell _{i}}^{i}$ be a shortest
nontrivial path (possibly of order two) between $V(C_{i})$ and $V(C_{0})$ in
$G$ for each $1\leq i\leq s,$ where $w_{\ell _{i}}^{i}\in V(C_{0})$. }If all
cycles $C_{1},\ldots ,C_{s}$ have length $\equiv 0\pmod 3$, then the
subgraphs $G_{2}=(C_{1}+P_{1}+C_{0})\cup C_{2}\cup \ldots \cup C_{s}$ and $%
G_{1}=G-G_{2}$ satisfy the conditions {and the result follows}. {Hence, we
assume that some }cycle $C_{i}\;(i\geq 1)$ has length $\equiv 2\pmod 3$. {%
Note that the }paths $P_{i}$'s {minus their end-vertices belonging to }$%
V(C_{0})$ are disjoint. {If some }$C_{i}$ has length $\equiv 0\pmod 3,$ say $%
C_{1}$, {then let }$L=C_{0}\cup (\cup _{i=2}^{s}C_{i})\cup (\cup
_{i=2}^{s}P_{i})$. By the induction hypothesis, $L$ has two subgraphs $%
L_{1},L_{2}$ such that $L_{1}$ has no cycle of $\mathcal{F}-\{C_{1}\}$ and
each component of $L_{2}$ belongs to $\mathcal{F}_{0}\cup \mathcal{F}%
_{0,2}\cup \mathcal{F}_{2,2}\cup \mathcal{E}$. Now $G_{2}=L_{2}\cup C_{1}$
and $G_{1}=G-G_{2}$ satisfy the conditions and the result holds. Hence we
can assume that all cycles $C_{1},\ldots ,C_{s}$ have length $\equiv 2%
\pmod 3
$. Let {$L=C_{0}\cup (\cup _{i=1}^{s}C_{i})\cup (\cup _{i=1}^{s}P_{i})$}. {%
Let $x_{i_{1}}^{0},\ldots ,x_{i_{t}}^{0}$ be the vertices of $C_{0}$ with
degree at least three and assume, without loss of generality, that $%
i_{1}<i_{2}<\ldots <i_{t}$.} Consider the following situations.

\smallskip \textbf{Subcase 2.1.} $t=2$.\newline
{If $\deg (x_{i_{1}}^{0}),\deg (x_{i_{2}}^{0})\geq 4$, then let $G_{2}$ be
the graph obtained from $L$ by deleting all vertices of $V(C_{0})-%
\{x_{i_{1}}^{0},x_{i_{2}}^{0}\}$. Otherwise, let $G_{2}$ be the graph
obtained from $L$ by deleting either the edge }$x_{i_{1}}^{0}x_{i_{2}}^{0}$ {%
(if any) or all the vertices $x_{i_{1}+1},\ldots ,x_{i_{2}-1}$. Then the
subgraphs $G_{2}$ and $G_{1}=G-G_{2}$ satisfies the conditions and the
result holds. }

\smallskip \textbf{Subcase 2.2.} $t=3$.\newline
{If $\deg (x_{i_{1}}^{0}),\deg (x_{i_{2}}^{0}),\deg (x_{i_{3}}^{0})\geq 4$,
then let $G_{2}$ be the graph obtained from $L$ by deleting all vertices of $%
V(C_{0})-\{x_{i_{1}}^{0},x_{i_{2}}^{0},x_{i_{3}}^{0}\}$. If $\deg
(x_{i_{1}}^{0})=\deg (x_{i_{2}}^{0})=\deg (x_{i_{3}}^{0})=3$, then let $%
G_{2} $ be the graph obtained from $L$ by deleting either the edge }$%
x_{i_{2}}^{0}x_{i_{3}}^{0}$ (if any) or all the vertices {of $%
\{x_{i_{2}+1}^{0},x_{i_{2}+2}^{0},\ldots ,x_{i_{3}-1}^{0}\}$. If, without
loss of generality, $\deg (x_{i_{1}}^{0})=3$ and $\deg (x_{i_{3}}^{0})\geq 4$%
. Let $G_{2}$ be the graph obtained from $L$ by deleting either the edge $%
x_{i_{2}}^{0}x_{i_{3}}^{0}$ (if any) or {all vertices between of }$%
x_{i_{2}}^{0}$ and $x_{i_{3}}^{0}$ {as well all vertices between $%
x_{i_{3}}^{0}$ and $x_{i_{1}}^{0} $ {starting from }$x_{i_{3}+1}^{0}.$ {In
either situation, the }subgraphs $G_{2}$ and $G_{1}=G-G_{2}$ satisfies the
conditions {and the result follows}. }}

\smallskip \textbf{Subcase 2.3.} $t\geq 4$.\newline
{If $\deg (x_{i_{1}}^{0}),\deg (x_{i_{2}}^{0}),\deg (x_{i_{3}}^{0})\geq 4$,
then let $G_{2}$ be the graph obtained from $L$ by deleting all vertices of $%
V(C_{0})-\{x_{i_{1}}^{0},x_{i_{2}}^{0},\ldots ,x_{i_{t}}^{0}\}$. If $\deg
(x_{i_{1}}^{0})=\deg (x_{i_{2}}^{0})=\ldots =\deg (x_{i_{t}}^{0})=3$, then
let $G_{2}$ be the graph obtained from $L$ by deleting all vertices of $%
\bigcup_{j=1}^{\lfloor t/2\rfloor
}\{x_{i_{2j}+1}^{0},x_{i_{2j}+2}^{0},\ldots ,x_{i_{2j+1}-1}^{0}\}$. Assume
without loss of generality that $\deg (x_{i_{1}}^{0})=3$ and $\deg
(x_{i_{t}}^{0})\geq 4$. Let $L^{1}$ be the component of $L-%
\{x_{i_{1}-1}^{0}x_{i_{1}}^{0},x_{i_{2}}^{0}x_{i_{2}+1}^{0}\}$ containing $%
x_{i_{1}}^{0}$, and let $L^{2}$ be the component of $L-%
\{x_{i_{3}-1}^{0}x_{i_{3}}^{0}\}$ containing $x_{i_{3}}^{0}$ if $\deg
(x_{i_{3}}^{0})\geq 4$, and be the component of $L-%
\{x_{i_{3}}^{0}x_{i_{3}-1}^{0},x_{i_{4}}^{0}x_{i_{4}+1}^{0}\}$ containing $%
x_{i_{3}}^{0}$ if $\deg (x_{i_{3}}^{0})=3$. Repeating this process we obtain
a sequence $L^{1},\ldots ,L^{p}$ of subgraphs $L$ which contains all cycles
of $L$ but $C_{0}$. Now the subgraphs $G_{2}=\cup _{i=1}^{p}L^{i}$ and $%
G_{1}=G-G_{2}$ satisfies the conditions and the result follows.}$\hfill \Box
$

\begin{lemma}
\label{mainn}Let $G$ be a connected graph with at least two disjoint cycles
of length $\equiv 0,2\pmod 3$, and let $\mathcal{F}$ be the family of all
cycles of $G$ with length $\equiv 0\;\mathrm{or}\;2\pmod 3$. Then there
exists a maximal subfamily $\mathcal{T}$ of pairwise disjoint cycles of $%
\mathcal{F}$ with $|\mathcal{T}|\geq 2$ and two disjoint subgraphs $G_{1}$
(possibly null) and $G_{2}$ of $G$ such that $V(G)=V(G_{1})\cup V(G_{2})$, $%
G_{1}$ has no cycle of $\mathcal{F}$ and each component of $G_{2}$ belongs
to $\mathcal{F}_{0}\cup \mathcal{F}_{0,2}\cup \mathcal{F}_{2,2}\cup \mathcal{%
E}$.
\end{lemma}

\noindent \textbf{Proof.} {By Lemma \ref{Lam3}, for any maximal subfamily }$%
\mathcal{T}$ of pairwise disjoint cycles of $\mathcal{F}$ with $|\mathcal{T}%
|\geq 2$, $G$ has two disjoint subgraphs $G_{1}^{\mathcal{T}}$ {and }$G_{2}^{%
\mathcal{T}}$ such that $V(G)=V(G_{1}^{\mathcal{T}})\cup V(G_{2}^{\mathcal{T}%
})$, $G_{1}^{\mathcal{T}}$ has no cycle of $\mathcal{T}$ {and }each
component of $G_{2}^{\mathcal{T}}$ is in $\mathcal{F}_{0}\cup \mathcal{F}%
_{0,2}\cup \mathcal{F}_{2,2}\cup \mathcal{E}. $ {Now, let }$c_{\mathcal{T}}$
{denote the }number of cycles of $G_{2}^{\mathcal{T}}$, {and let }$s_{%
\mathcal{T}}$ be the sum of the lengths of paths between two cycles in the
components of $G_{2}^{\mathcal{T}}$ {that belong to }$\mathcal{F}_{0,2}\cup
\mathcal{F}_{2,2}\cup (\cup _{r,s\geq 0;r+s=2}B_{r,s}).$ {\ Moreover, let }
\begin{equation*}
c_{\mathcal{F}}=\max \{c_{\mathcal{T}}\mid \mathcal{T}\;\mathrm{%
is\;a\;maximal\;subfamily\;of\;pairwise\;disjoint\;cycles\;of}\;\mathcal{F}\;%
\mathrm{with}\;|\mathcal{T}|\geq 2\}.
\end{equation*}%
Choose a triple $(\mathcal{T},G_{1}^{\mathcal{T}},G_{2}^{\mathcal{T}})$ such
that: (i) $c_{\mathcal{F}}=c_{\mathcal{T}}$; (ii) subject to (i): $s_{%
\mathcal{T}}$ is maximized. {Notice that }$G_{2}^{\mathcal{T}}$ {may not
contain all cycles of }$\mathcal{T}$. We claim that {the }two disjoint
subgraphs $G_{1}^{\mathcal{T}}$ and $G_{2}^{\mathcal{T}}$ {chosen in this
way yield the desired result. }

{It is clear that it suffices to show }that $G_{1}^{\mathcal{T}}$ has no
cycle of $\mathcal{F}$. {Hence, suppose }to the contrary that $G_{1}^{%
\mathcal{T}}$ contains at least one cycle of $\mathcal{F}$. Let $G_{2}^{1}$\
be obtained from $G_{2}^{\mathcal{T}}$ by adding { a maximum set
of pairwise of cycles of $\mathcal{F}$ } with length $\equiv 0\pmod 3$ {\
belonging to }$G_{1}^{\mathcal{T}}$ and let $G_{1}^{1}=G-G_{2}^{1}$. Note
that if $G_{1}^{\mathcal{T}}$ {contains no }cycle of $\mathcal{F} $ with
length $\equiv 0\pmod 3$, then $G_{2}^{1}=G_{2}^{\mathcal{T}}$. {Now, l}et $%
\mathcal{T}_{1}^{1}$ be the family of all cycles of $G_{2}^{1}$ {that belong
}to $\mathcal{F}$ and {et }$\mathcal{T}^{1}$ be a maximal subfamily of $%
\mathcal{F}$ such that $\mathcal{T}_{1}^{1}\subseteq \mathcal{T}^{1}$. If $%
G_{1}^{1}$ does not contain any cycle of $\mathcal{F}$, then the family $%
\mathcal{T}^{1}$ and the subgraphs $G_{1}^{1}$ and $G_{2}^{1}$ satisfy the
conditions {which }leads to a contradiction because of $c_{\mathcal{T}%
^{1}}>c_{\mathcal{T}}$. {Hence we }assume that $G_{1}^{1}$ contains at least
one cycle of $\mathcal{F}$.

Let first $H^{1},\ldots H^{r}$ be the components of $G_{1}^{1}$ which
contains at least two disjoint cycles of $\mathcal{F}$ (if any), and let $%
\mathcal{F}^{i}$ be a maximal subfamily of pairwise disjoint cycles of $%
\mathcal{F}$ {that are }in $H^{i}$ with $|\mathcal{F}^{i}|\geq 2$, for each $%
1\leq i\leq r$. By Lemma \ref{Lam3}, $H^{i}$ has two subgraphs $%
H_{1}^{i},H_{2}^{i}$ such that $H_{1}^{i}$ has no cycles of $\mathcal{F}^{i}$
and each component of $H_{2}^{i}$ is in $\mathcal{F}_{0}\cup \mathcal{F}%
_{0,2}\cup \mathcal{F}_{2,2}\cup \mathcal{E}$. Let $G_{2}^{2}=G_{2}^{1}\cup
(\cup _{i=1}^{r}H_{2}^{i})$, if $r\geq 1$ and $G_{1}^{2}=G-G_{2}^{2}$. Let $%
\mathcal{T}_{2}^{1}$ be the family of all cycles of $G_{2}^{2}$ {\ that
belong }to $\mathcal{F}$ and let $\mathcal{T}^{2}$ be a maximal subfamily of
$\mathcal{F}$ such that $\mathcal{T}_{1}^{2}\subseteq \mathcal{T}^{2}$. If $%
G_{1}^{2}$ does not contain any cycle of $\mathcal{F}$, then the family $%
\mathcal{T}^{2}$ and the subgraphs $G_{1}^{2}$ and $G_{2}^{2}$ satisfy the
conditions {which }leads to a contradiction because $c_{\mathcal{T}^{2}}>c_{%
\mathcal{T}}$. {Hence, we }assume that $G_{1}^{2}$ contains at least one
cycle of $\mathcal{F}$. If $G_{1}^{2}$ has a component with at least two
disjoint cycles of $\mathcal{F}$, then we proceed as above.

Henceforth, we {can }assume that {each }component of $G_{1}^{2}$ has at most
one cycle of $\mathcal{F}$. Let $C_{0}$ be a cycle of $G_{1}^{2}$ {belonging
}to $\mathcal{F}$. Clearly, $C_{0}=x_{1}^{0}x_{2}^{0}\ldots
x_{m_{0}}^{0}x_{1}^{0}$ is connected to a component of $G_{2}^{2}$ by {some }%
path (possibly an edge). Let $P=(x_{1}^{0}=)v_{0}v_{1}\ldots v_{t}$ be a
shortest path between $V(C_{0})$ and $V(G_{2}^{2})$. Then $v_{t}$ belongs to
{a component of }$G_{2}^{2}$, {say }$H_{1}.$ If $H_{1}\in \mathcal{F}_{0}$,
then let $H_{1}^{\prime }=H_{1}+P+C_{0}$ and let $%
G_{2}^{3}=(G_{2}^{2}-H_{1})\cup H_{1}^{\prime }$. Hence, assume that $H_{1}$
contains at least two cycles. We distinguish the following cases.

\smallskip \textbf{Case 1.} $H_{1}\in \mathcal{F}_{2,2} $.%
\newline
Suppose $H_{1}$ is obtained from two cycles $C_{1}=x_{1}^{1}\ldots
x_{m_{1}}^{1}x_{1}^{1}$ and $C_{2}=x_{1}^{2}\ldots x_{m_{2}}^{2}x_{1}^{2}$
by adding an edge $x_{1}^{1}x_{1}^{2}$. We further assume, without loss of
generality, that $v_{t}=x_{j}^{2}\in V(C_{2})$. {Let }$H_{1}^{\prime }$ be
the graph obtained from $C_{0}\cup C_{1}$ to which we add the path $%
Px_{j-1}^{2}\ldots x_{2}^{2}Q$, {in other words, $H_{1}^{\prime }$ {is
obtained from }$C_{0}+P+H_{1}$ by removing vertices $%
x_{j+1}^{2},...,x_{m_{2}}^{2}.$ Note that }$H_{1}^{\prime }$ belongs to
either $\mathcal{B}_{0,2}\cup \mathcal{B}_{1,1}.$ Now let $%
G_{2}^{3}=(G_{2}^{2}-H_{1})\cup H_{1}^{\prime} $ which we will discuss
further below.

{ \smallskip \textbf{Case 2.} $H_{1}\in \mathcal{F}_{0,2} $.%
\newline
Suppose $H_{1}$ is obtained from two cycles $C_{1}=x_{1}^{1}\ldots
x_{m_{1}}^{1}x_{1}^{1}$ and $C_{2}=x_{1}^{2}\ldots x_{m_{2}}^{2}x_{1}^{2}$
by adding a path $Q=(x_{1}^{1}=)z_{0}z_{1}\ldots z_{k}(=x_{1}^{2})$. Suppose
without loss of generality that $m_1\equiv 0\pmod 3$ and $m_2\equiv \pmod 3$%
. If $v_t=z_j$ for some $j$, then let $H_1^{\prime}$ be obtained from $%
C_{0}\cup C_{2}$ to which we add the path $Pz_j\ldots z_{k}$, and let $%
G_{2}^{3}=(G_{2}^{2}-H_{1})\cup (H_{1}^{\prime} \cup C_1)$. Suppose that $%
v_t\in V(C_1)\cup V(C_2)$. {We further assume, without loss of generality,
that }$v_{t}=x_{j}^{2}\in V(C_{2})$. {Let }$H_{1}^{\prime }$ {be the graph
obtained from }$C_{0}\cup C_{1}$ {to which we add the path }$%
Px_{j-1}^{2}\ldots x_{2}^{2}Q$, {in other words, }$H_{1}^{\prime }$ {is
obtained from }$C_{0}+P+H_{1}$ {by removing vertices }$%
x_{j+1}^{2},...,x_{m_{2}}^{2}.$ {Note that }$H_{1}^{\prime }$ belongs to
either $\mathcal{F}_{0,2}$ or $\mathcal{B}_{0,2}\cup \mathcal{B}_{1,1}.$ Now
let $G_{2}^{3}=(G_{2}^{2}-H_{1})\cup H_{1}^{\prime } $ which we will discuss
further below. }

\smallskip \textbf{Case 3.} $H_{1}\in \mathcal{B}_{r,s}$ where $r+s\geq 2$.%
\newline
Let $z^{\ast }$ {be the special vertex of }$H_{1}.$ {If }$v_{t}=z^{\ast }$,
then $H_{1}^{\prime }=H_{1}+P+C_{0}$ is a subgraph belonging to $\mathcal{E}$%
. In this case, let $G_{2}^{3}=(G_{2}^{2}-H_{1})\cup H_{1}^{\prime }$ which
we will discuss further below. Hence we assume that $v_{t}\neq z^{\ast }.$
First let $r+s=2.$ Then $H_{1}$ is obtained from two cycles $%
C_{1}=x_{1}^{1}\ldots x_{m_{1}}^{1}x_{1}^{1}$ and $C_{2}=x_{1}^{2}\ldots
x_{m_{2}}^{2}x_{1}^{2}$ by adding a path $Q=(x_{1}^{1}=)z_{0}z_{1}\ldots
z_{t}(=x_{1}^{2}),$where $t\geq 2.$ If $v_{t}\in \{z_{1},\ldots ,z_{t-1}\}$,
then let $H_{1}^{\prime }=H_{1}+C_{0}+P$ and $G_{2}^{3}=(G_{2}^{2}-H_{1})%
\cup H_{1}^{\prime }$ . {Now, suppose, without loss of generality, that }$%
v_{t}=x_{j}^{2}\in V(C_{2}). $ {Let }$H_{1}$ be obtained {from }$C_{0}\cup
C_{1}$ {to which we add the path }$Px_{j-1}^{2}\ldots x_{2}^{2}Q.$ S{et }$%
G_{2}^{3}=(G_{2}^{2}-H_{1})\cup H_{1}^{\prime }$ {which we will discuss
further below}.

Now let $r+s\geq 3.$ {Assume that }$v_{t}$ {belongs to one of the }${s+r}$
cycles of $H_{1},$ {say} $C^{\prime }.$ Let $H_{1}^{\prime }=C^{\prime
}+C_{0}+P$ and $H_{1}^{\prime \prime }$ be the graph obtained from $H_{1}$
by deleting the vertices of $V(C^{\prime })$ and the path (if any) joining $%
z^{\ast }$ to $V(C^{\prime })$ in $H_{1}$. Note that $H_{1}^{\prime }$
belongs to either $\mathcal{F}_{0,2}$ or $\mathcal{B}_{0,2}\cup \mathcal{B}%
_{1,1}.$ Now let $G_{2}^{3}=(G_{2}^{2}-H_{1})\cup H_{1}^{\prime }$ . Assume
that $v_{t}$ belongs to a path on a tailed cycle $C_{m,\ell }$ of $H_{1},$
and let $P^{\prime }$ be the subpath between $v_{t}$ and the cycle $%
C^{\prime }$ of $C_{m,\ell }.$ {Let }$H_{1}^{\prime }=C^{\prime
}+C_{0}+P+P^{\prime }$ and $H_{1}^{\prime \prime }$ be the graph obtained
from $H_{1}$ by deleting {the }vertices of $V(C_{m,\ell }).$ {Note that }$%
H_{1}^{\prime }$ {belongs to either }$\mathcal{F}_{0,2}$ or $\mathcal{B}%
_{0,2}.$ Suppose $G_{2}^{3}=(G_{2}^{2}-H_{1})\cup H_{1}^{\prime }\cup
H_{1}^{\prime \prime }$. Obviously either the number of cycles of $G_{2}^{3}$
is greater than the number of cycles of $G_{2}^{\mathcal{T}}$ or the sum of
lengths of paths between two cycles of $G_{2}^{3}$ that belong to $\mathcal{F%
}_{0,2}\cup \mathcal{F}_{2,2}\cup (\cup _{r,s\geq 0;r+s=2}B_{r,s})$ is
greater than the corresponding sum of $G_{2}^{\mathcal{T}}$. Let $%
G_{1}^{3}=G-G_{2}^{3}.$ Let $\mathcal{T}_{3}^{1}$ be the family of all
cycles of $G_{2}^{3}$ which belongs to $\mathcal{F}$ and let $\mathcal{T}%
^{3} $ be a maximal subfamily of $\mathcal{F}$ such that $\mathcal{T}%
_{1}^{3}\subseteq \mathcal{T}^{3}$. If $G_{1}^{3}$ does not contain any
cycle of $\mathcal{F}$, then the family $\mathcal{T}^{3}$ and the subgraphs $%
G_{1}^{3}$ and $G_{2}^{3}$ satisfy the conditions which leads to a
contradiction because {of either }$c_{\mathcal{T}^{3}}>c_{\mathcal{T}}$ or $%
s_{\mathcal{T}^{3}}>s_{\mathcal{T}}$. Hence we assume that $G_{1}^{3}$
contains at least one cycle of $\mathcal{F}$. We repeat the above precess.
Since $G$ is finite, this process will stop and we obtain a maximal
subfamily $\mathcal{T}^{\prime }$ of pairwise disjoint cycles of $\mathcal{F}
$ with $|\mathcal{T}^{\prime }|\geq 2 $ and two disjoint subgraphs $G_{1}$
(possibly null), $G_{2}$ of $G$ such that $V(G)=V(G_{1})\cup V(G_{2})$, $%
G_{1}$ has no cycle of $\mathcal{F}$ and each component of $G_{2}$ belongs
to $\mathcal{F}_{0}\cup \mathcal{F}_{0,2}\cup \mathcal{F}_{2,2}\cup \mathcal{%
E}$. $\hfill \Box $

\begin{lemma}
\label{mainnn}\emph{Let $G$ be a connected graph with at least two disjoint
cycles of length $\equiv 0,2\pmod 3$, and let $\mathcal{F}$ be the family of
all cycles of $G$ with length $\equiv 0\;\mathrm{or}\;2\pmod 3$. Then there
exists a maximal subfamily $\mathcal{T}$ of pairwise disjoint cycles of $%
\mathcal{F}$ with $|\mathcal{T}|\geq 2$ and two disjoint subgraphs $G_{1}$
(possibly null), $G_{2}$ of $G$ such that $V(G)=V(G_{1})\cup V(G_{2})$, $%
G_{1}$ has no cycle of $\mathcal{F}$ and each component of $G_{2}$ belongs
to $\mathcal{F}_{0}\cup \mathcal{F}_{0,2}\cup \mathcal{F}_{2,2}\cup \mathcal{%
F}_{3}\cup \mathcal{E}$.}
\end{lemma}

\noindent \textbf{Proof.} Let $(\mathcal{T},G_{1},G_{2})$ {be the triple
satisfying the conditions of Lemma \ref{mainn}. Hence }$G_{1}$ has no cycle
of $\mathcal{F}$ {and each component of }$G_{2}$ {belongs to }$\mathcal{F}%
_{0}\cup \mathcal{F}_{0,2}\cup \mathcal{F}_{2,2}\cup \mathcal{E}$. {If there
are no two components }$H_{1}\in \mathcal{F}_{0,2},$ $H_{2}\in \mathcal{F}%
_{2,2}$ of $G_{2}$ {joined by a path }$P$ {in }$G$ {with all its vertices,
except the end-vertices, belong to }$V(G_{1})$, then $G_{1}$ and $G_{2}$ are
{the }desired subgraphs. {Hence we assume }that there {are two components }$%
H_{1}\in \mathcal{F}_{0,2}$ and $H_{2}\in \mathcal{F}_{2,2}$ of $G_{2}$ {%
joined by a path }$P$ {in }$G$ {with all its vertices, except the
end-vertices, belong to }$V(G_{1}).$ Let $G_{2}^{\prime }$ be {\ the graph }%
obtained from $G_{2}$ by adding the path $P$ and let $G_{1}^{\prime
}=G-G_{2}^{\prime }$. Clearly $G_{1}^{\prime }$ and $G_{2}^{\prime }$
satisfy the conditions {and the result follows}. We can repeat this process {%
until we get two }subgraphs $G_{1}^{\ast }$ (possibly null) and $G_{2}^{\ast
}$ such that $V(G)=V(G_{1}^{\ast })\cup V(G_{2}^{\ast })$, $G_{1}^{\ast }$
has no cycle of $\mathcal{F}$, each component of $G_{2}^{\ast }$ is in $%
\mathcal{F}_{0}\cup \mathcal{F}_{0,2}\cup \mathcal{F}_{2,2}\cup \mathcal{F}%
_{3}\cup \mathcal{E}$ and {such that no path in }$G$ {like to the one
described above joins two components }$H^{\prime }\in \mathcal{F}_{0,2}$ and
$H^{\prime \prime }\in \mathcal{F}_{2,2}$ of $G_{2}^{\ast }.$ $\hfill \Box $

\bigskip

{From now on, }a graph in $\mathcal{(}\cup _{r+s\geq 2;s\leq 2}\mathcal{B}%
_{r,s})\cup \mathcal{F}_{3}\cup \mathcal{F}_{0,2}\cup \mathcal{F}_{0}$ {will
be called }\emph{strong}. {Also, the }special vertex of each graph in $%
\mathcal{B}_{r,s}$ {will be }called a \textit{strong vertex}.

\begin{lemma}
\label{strong}\emph{Let $k\geq 1$ be an integer and let $G$ be a graph of
order $n$ and minimum degree $\delta \geq 2$, which does not contain any
induced $\{C_{5},C_{8},\ldots ,C_{3k+2}\}$-cycles. If $G$ {is strong, }then $%
G$ has a 3-tuple $\overrightarrow{f}$ of {RDFs} such that $\omega (%
\overrightarrow{f})\leq \frac{(4k+8)3n}{6k+11}$ and all vertices of $G$ are $%
\overrightarrow{f}$-strong. }
\end{lemma}

\noindent \textbf{Proof.} Let $G\in (\cup _{r+s\geq 2;s\leq 2}\mathcal{B}%
_{r,s})\cup \mathcal{F}_{3}\cup \mathcal{F}_{0,2}\cup \mathcal{F}_{0}.$
Assume first that $G\in \mathcal{F}_{0}$. Then $\gamma _{R}(G)=\frac{2n}{3}<%
\frac{(4k+8)n}{6k+11}$. Let $G=x_{1}x_{2}\ldots x_{3t}x_{1}$ and define for $%
j\in \{1,2,3\}$ the functions $f_{j}$ on $V(G)$ as follows: $%
f_{j}(x_{3i+j})=2$ for $0\leq i\leq t-1$ and $f_{j}(x)=0$ otherwise. Clearly
$f_{j}$ is an $\gamma _{R}(G)$-function for each $j\in \{0,1,2\}$ and the
triple $\overrightarrow{f}=(f_{0},f_{1},f_{2})$ satisfies the desired result.

Assume now that $G\in \mathcal{F}_{0,2}$. {Since }$G$ has no induced $%
\{C_{5},C_{8},\ldots ,C_{3k+2}\}$-cycles, {we deduce that cycle of lentgh }$%
\equiv 2\pmod 3$ in $\mathcal{F}_{0,2}${\ has order at least }$3k+5$, and
thus $G$ has order at least $3k+8.$ {Now b}y Lemma \ref{MainLem}, $G$ has a
3-tuple $\overrightarrow{f}$ of RDFs such that $\omega (\overrightarrow{f}%
)\leq 2n+1$ and all vertices of $G$ are $\overrightarrow{f}$-strong. {A
simple calculation shows that }$\omega (\overrightarrow{f})\leq \frac{%
(4k+8)3n}{6k+11}$.

Next assume that\textbf{\ }$G\in \mathcal{F}_{3}$. {By definition, }$G$ is
obtained from a graph $G_{1}\in \mathcal{F}_{0,2}$ and a graph in $G_{2}\in
\mathcal{F}_{2,2}$ by adding either an edge $vw$ {or a path }$Q$ {joining a
vertex of }$G_{1}$ to a vertex of $G_{2}$ {so that all vertices of }$Q${\
become of degree two in }$G.$ Let $G_{1}$ be obtained from two cycles $%
C_{1}=x_{1}^{1}x_{2}^{1}\ldots x_{n_{1}}^{1}x_{1}^{1}\in \mathcal{F}_{0}$
and $C_{2}=x_{1}^{2}x_{2}^{2}\ldots x_{n_{2}}^{2}x_{1}^{2}\in \mathcal{F}%
_{2} $ by adding {\ either the edge }$x_{1}^{1}x_{1}^{2}$ {or a }path $P$
between $x_{1}^{1}$ and $x_{1}^{2}.$ {By Lemma \ref{MainLem} (items 8,9,10),
$G_{1}$ has a a 3-tuple $\overrightarrow{f}=(f_{1},f_{2},f_{3})$ of RDFss
such that $\omega (\overrightarrow{f},G_{1})\leq 2n(G_{1})+1$ and all
vertices of $G_{1} $ are $\overrightarrow{f}$-strong. Moreover, }let $G_{2}$
be obtained from two cycles $C_{3}=x_{1}^{3}x_{2}^{3}\ldots
x_{n_{1}}^{3}x_{1}^{3}\in \mathcal{F}_{2}$ and $C_{4}=x_{1}^{4}x_{2}^{4}%
\ldots x_{n_{4}}^{4}x_{1}^{4}\in \mathcal{F}_{2}$ by adding the edge $%
x_{1}^{3}x_{1}^{4}$. Without loss of generality, we assume that the added
edge $uv$ or the path $Q$ is between $V(C_{3})$ and $V(G_{1})$. {By
sequentially applying Lemmas \ref{tailedcycle-3p+1} (items 3,4) (once on }$%
uv $ or $Q$ and $C_{3},$ {and then on the resulting graph with }$C_{4}),$ {$%
\overrightarrow{f}$ can be extended to a triple $\overrightarrow{g}$ of RDFs
of $G$ such that $\omega (\overrightarrow{g},G_{2})\leq 2n(G_{2})+2$ and
each vertex of $G_{2}$ is $\overrightarrow{g}$-strong. Since }$G$ has no
induced $\{C_{5},C_{8},\ldots ,C_{3k+2}\}$-cycles, {we deduce that order
each cycle of lentgh }$\equiv 2\pmod 3$ in $G$ is at least $3k+5.$ {Using
the fact that }$G$ has three cycles {of length }$\equiv 2\pmod 3$ and one
cycle {of lentgh }$\equiv 0\pmod 3,$ we have {$n(G)\geq 9k+18.$ Therefore $%
\omega (\overrightarrow{g})\leq 2n(G)+3\leq \frac{(4k+8)3n(G)}{6k+11}.$ }

Using a similar argument we can show that for any graph $G\in \mathcal{\cup }%
_{r+s\geq 2;s\leq 2}B_{r,s}$ the result is {also }true. $\hfill \Box $

\begin{lemma}
\label{non-strong1}\emph{Let $k\geq 1$ be an integer and let $G$ be a graph
of order $n$, minimum degree $\delta \geq 2$, which does not contain any
induced $\{C_{5},C_{8},\ldots ,C_{3k+2}\}$-cycles. If $G\in \mathcal{B}%
_{r,s} $ with $s\geq 3$, then $G$ has a 3-tuple $\overrightarrow{f}$ of RDFs
such that $\omega (\overrightarrow{f})\leq \frac{(4k+8)3n}{6k+11}$ and the
special vertex as well as all vertices on tailed cycles of $G$ are $%
\overrightarrow{f}$-strong. }
\end{lemma}

\noindent \textbf{Proof.} Suppose $G$ be obtained from $r\geq 0$ graphs $%
C_{n_{1},\ell _{1}},\ldots ,C_{n_{r},\ell _{r}}$ and $s\geq 3$ cycles $%
C_{m_{1}},\ldots ,C_{m_{s}}$, where $n_{i}\equiv 2\pmod 3$ and $m_{j}\equiv 2%
\pmod 3$ for each $i,j$, by adding a new vertex $z$ (special vertex)
attached to {endvertices of the }$C_{n_{i},\ell _{i}}$'s and to one vertex
of each cycle $C_{m_{j}}$. {We first note that each of the }$r+s\geq 3$ {%
cycles has order at least }$3k+5$, {and thus each tailed cycle contains at
least }${(3k+5)+1}$ vertices. {Hence }${n(G)\geq (3k+5)(s+r)+r+1.}$ {\ Now,
if }$r=0$, then the result follows from Lemmas \ref{MainLem}-(11) {and the
previous fact}. {Hence assume }that $r\geq 1$. Let $H$ be obtained from $G$
by deleting all vertices of $C_{n_{i},\ell _{i}}$'s. By Lemma \ref{MainLem}
(item 11), $H$ has a triple $\overrightarrow{f}$ such that $\omega (%
\overrightarrow{f})\leq 2n(H)-s+4$ and $z$ is $\overrightarrow{f}$-strong. {%
Since }${n(H)\geq (3k+5)s+1,}$ we deduce that $\omega (\overrightarrow{f}%
)\leq \frac{(4k+8)3n(H)}{6k+11}$. {Now, by applying repeatedly }Lemma \ref%
{tailedcycle-3p+1}-(4) {on }$C_{n_{1},\ell _{1}},\ldots ,C_{n_{r},\ell _{r}}$%
, we can extend $\overrightarrow{f}$ to a triple $\overrightarrow{g}$ of $G$
such that $\omega (g,\cup _{i=1}^{r}C_{n_{1},\ell _{1}})\leq
\sum_{i=1}^{r}(2n(C_{n_{1},\ell _{1}})+1)$ and all newly added vertices are $%
\overrightarrow{g}$-strong. {Therefore, }$\omega (\overrightarrow{g})\leq
2n+r+4-s.$ {Now by the previous fact on the order and the calculation, we
can see that }$2n+r+4-s\leq \frac{(4k+8)3n}{6k+11}$, {which proves the
result.}$\hfill \Box ${\ }

\begin{lemma}
\label{non-strong2} Let $k\geq 1$ be an integer and let $G\in \mathcal{F}%
_{2,2}$ be a graph of order $n$, minimum degree $\delta \geq 2$, which does
not contain any induced $\{C_{5},C_{8},\ldots ,C_{3k+2}\}$-cycles. Then

\begin{enumerate}
\item $G$ has a $3$-tuple $f$ of RDFs such that $\overrightarrow{f}\leq \
2n(G) + 1 \leq \frac{(4k+8)3n(G)}{(6k+11)}.$

\item If H is a graph obtained from $G$ and a cycle $%
C_{3p+1}=x_{1}...x_{3p+1}x_1$ by adding an edge between them, then H has a
3-tuple $f$ of RDFs such that $\omega(\overrightarrow{f}) \leq \frac{%
(4k+8)3n(G)}{(6k+11)}$ and all vertices of H but $x_{3p+1}$ are $%
\overrightarrow{f}$-strong.

\item If $H$ is a graph obtained from $G$ and a tailed cycle $C_{3p+1,\ell }$
with vertex set $x_1,\ldots ,x_{3p+1}$, $y_{1},\ldots ,y_{\ell }$ by joining
$y_{\ell }$ to a vertex of G, then H has a 3-tuple $f $ of RDFs such that $%
\omega(\overrightarrow{f}) \leq \frac{\emph{(4k+8)3n(G)}}{(6k+11)}$ and all
vertices of $H$ but $x_{3p+1}$ are $\overrightarrow{f}$-strong.
\end{enumerate}
\end{lemma}

\noindent \textbf{Proof.} (1) is easy to show and so we prove only (2) and
(3). Let $G\in \mathcal{F}_{2,2}$ be formed from two cycles $C_{1}$ and $%
C_{2}$ by adding an edge between them, and let $H$ be obtained from $G$ and
the cycle $C_{3p+1}$ (resp. tailed cycle $C_{3p+1,\ell }$) by adding an edge
$xy$ (resp. $xy_{\ell }$), where without loss of generality $x\in V(C_{2})$.
Let $K$ be the graph obtained from $H$ by deleting all vertices of $V(C_{1})$%
. By Lemma \ref{MainLem} (items 2,3 and 4), $K$ has a 3-tuple $%
\overrightarrow{g}$ of RDFs of $K$ such that $\omega (\overrightarrow{g}%
)\leq 2n(K)+1$ and all vertices of $K$ except $x_{3p+1}$ are $%
\overrightarrow{g}$-strong. Now by Lemma \ref{tailedcycle-3p+1}, we can
extend $\overrightarrow{g}$ to a 3-tuple $\overrightarrow{f}$ of RDFs of $H$
such that $\omega (\overrightarrow{f})\leq 2n(H)+2$ and all vertices of $H$
except $x_{3p+1}$ are $\overrightarrow{g}$-strong. By assumption we have $%
n(H)\geq 6k+14$ and {thus one can check that }$\omega (\overrightarrow{f}%
)\leq \frac{(4k+8)3n(H)}{6k+11}$. $\hfill\Box $

\section{\protect\smallskip Proof of Conjecture \protect\ref{conj}}

Now we are ready to state our main result.{\ }

\begin{theorem}
\label{Theorem} \emph{Let $G$ be a graph of order $n\ge 6k + 9$, minimum
degree $\delta\ge 2$, which does not contain any induced $%
\{C_5,C_8,\ldots,C_{3k+2}\}$-cycles. Then $\gamma_R(G)\le \frac{(4k+8)n}{%
6k+11}.$}
\end{theorem}

\noindent \textbf{Proof.} Let $\mathcal{F}$ be the family of all cycles of $%
G $ with length $\equiv 0\;\mathrm{or}\;2\pmod 3$. If $|\mathcal{F}|=0$,
then the result follows from Theorem \ref{Th1} and if $|\mathcal{F}|\geq 1$
and $\mathcal{F}$ contains a cycle which intersect any cycle of $\mathcal{F}$%
, then the result follows from Theorems \ref{Th2} and \ref{Th3}. Henceforth,
we assume that each cycle of $\mathcal{F}$ belongs to a maximal subfamily $%
\mathcal{T}$ of pairwise disjoint cycles of $\mathcal{F}$ with $|\mathcal{T}%
|\geq 2$. Let $(G_{1}^{1},G_{2}^{1}),\ldots ,(G_{1}^{m},G_{2}^{m})$ be all
pairs of subgraph such that $V(G)=V(G_{1}^{i})\cup V(G_{2}^{i})$, $G_{1}^{i}$
has no cycle of $\mathcal{F}$ and each component of $G_{2}^{i}$ belongs to $%
\mathcal{F}_{0}\cup \mathcal{F}_{0,2}\cup \mathcal{F}_{2,2}\cup \mathcal{F}%
_{3}\cup \mathcal{E}$. Let{\ }$s_{(G_{1}^{i},G_{2}^{i})}$ be the sum of the
lengths of paths between two cycles in the components of $G_{2}^{\mathcal{T}%
} $ {that belong to }$\mathcal{F}_{0,2}\cup \mathcal{F}_{2,2}\cup (\cup
_{r,s\geq 0;r+s=2}B_{r,s}).$ {Among all pairs }$(G_{1}^{i},G_{2}^{i}),$ let $%
(G_{1},G_{2})$ {be one chosen so that: }

\begin{description}
\item[(C$_{1}$)] the {number of }strong components of $G_{2}$ is maximized.

\item[(C$_{2}$)] subject to Condition (C$_{1}$): the number of cycles of $%
G_{2}$ belonging to $\mathcal{F}$ is maximized.

\item[(C$_{3}$)] subject to Conditions (C$_{1}$) and (C$_{2}$): the number
of components of $G_{2}$ in $\mathcal{F}_{2,2}$ is minimized.

\item[(C$_{4}$)] subject to Conditions (C$_{1}$), (C$_{2}$) and (C$_{3}$): $%
s_{(G_{1},G_{2})}$ is maximized.
\end{description}

We proceed with some {further claims that are needed for our proof. }

\smallskip \noindent \textbf{Claim 1.} Let $M$ be a component of $G_{2}$
such that $M\in \mathcal{F}_{2,2}$. Then there is no path $v_{0}v_{1}\ldots
v_{t+1}\;(t\geq 1)$ in $G$ such that $v_{0},v_{t+1}\in M$, $v_{1},\ldots
,v_{t}\in V(G_{1})$ and $v_{0}$ and $v_{t+1}$ belong to different cycles of $%
M$.\newline
\noindent \textbf{Proof of Claim 1.} Suppose, to the contrary, that there is
path $P=v_{0}v_{1}\ldots v_{t+1}\;(t\geq 1)$ in $G$ such that $%
v_{0},v_{t+1}\in M$, $v_{1},\ldots ,v_{t}\in V(G_{1})$ and $v_{0}$ and $%
v_{t+1}$ belong to different cycles of $M$. {Let }$e^{\ast }$ {be the edge
joining the two cycles of }$M${\ and let }$M^{\prime }$ be obtained from $M$
by deleting $e^{\ast }$ and adding path $P$. {Set }$G_{2}^{\prime
}=(G_{2}-M)\cup M^{\prime }$ and $G_{1}^{\prime }=G-G_{2}$. Clearly $%
V(G)=V(G_{1}^{\prime })\cup V(G_{2}^{\prime })$, $G_{1}$ has no cycle of $%
\mathcal{F}$ and each component of $G_{2}^{\prime }$ {is in }$\mathcal{F}%
_{0}\cup \mathcal{F}_{0,2}\cup \mathcal{F}_{2,2}\cup \mathcal{F}_{3}\cup
\mathcal{E}$. But $G_{2}^{\prime }$ has one more strong component than $%
G_{2},$ {contradicting our choice of }$(G_{1},G_{2}).$ $\ \ \ \ \ \ \ \ \ \
\ \hspace{3cm}\Box $ \ \ \

\bigskip\

\smallskip \noindent \textbf{Claim 2.} {For any two components }$M_{1}$ and $%
M_{2}$ of $G_{2}$ {belonging to }$\mathcal{F}_{2,2}$, {there }is no path $%
v_{0}v_{1}\ldots v_{t+1}\;(t\geq 1)$ in $G$ such that $v_{0}\in M_{1}$, $%
v_{t+1}\in M_{2}$ and $v_{1},\ldots ,v_{t}\in V(G_{1})$.\newline
\noindent \textbf{Proof of Claim 2.} Suppose, to the contrary, that {for }%
two components $M_{1},M_{2}$ of $G_{2}$ {\ belonging to }$\mathcal{F}_{2,2},$
there is a path $v_{0}v_{1}\ldots v_{t+1}\;(t\geq 1)$ in $G_{1}$ such that $%
v_{0}\in M_{1}$ and $v_{t+1}\in M_{2}$. Suppose {that } $M_{1}$ is obtained
from two cycles $C_{1}=u_{1}^{1}\ldots u_{m_{1}}^{1}u_{1}^{1}$ and $%
C_{2}=u_{1}^{2}\ldots u_{m_{2}}^{2}u_{1}^{2}$ by adding {the edge }$%
u_{1}^{1}u_{1}^{2}$, and let $M_{2}$ {be }obtained from two cycles $%
C_{3}=u_{1}^{3}\ldots u_{m_{3}}^{3}u_{1}^{3}$ and $C_{4}=u_{1}^{4}\ldots
u_{m_{4}}^{4}u_{1}^{4}$ by adding the edge $u_{1}^{3}u_{1}^{4}$. {Moreover,
assume, without loss of generality, that }$v_{0}=u_{j}^{2}\in V(C_{2})$
where $j\geq m_{2}/2$ ({by relabeling the vertices if necessary) }and $%
v_{t+1}=u_{b}^{3}\in V(C_{3})$ where $b\geq m_{3}/2$ ({by relabeling the
vertices if necessary)}. {Now, let }$M$ be the subgraph obtained from $C_{1}$
and $C_{4} $ by adding the path $u_{1}^{1}u_{1}^{2}u_{2}^{2}\ldots
u_{j}^{2}v_{1}v_{2}\ldots v_{t}u_{b}^{3}u_{b-1}^{3}\ldots u_{1}^{3}u_{1}^{4}$%
. {Set }$G_{2}^{\prime }=(G_{2}-(M_{1}\cup M_{2}))\cup M$ and $G_{1}^{\prime
}=G-G_{2}^{\prime }$. If $G_{1}^{\prime }$ has no cycle of $\mathcal{F}$,
then by considering the pair $(G_{1}^{\prime },G_{2}^{\prime }) $ we get {%
one more strong component in }$G_{2}^{\prime }$ {than in }$G_{2},$ {%
contradicting our choice of }$(G_{1},G_{2}).$ {Hence we }assume that $%
G_{1}^{\prime }$ has some cycles of $\mathcal{F}$.

First let $G_{1}^{\prime }$ has exactly one cycle $C$ of $\mathcal{F}$. If $%
C $ has length $\equiv 0\pmod 3$, then as above we get a contradiction by
considering the subgraphs $G_{2}^{\prime \prime }=G_{2}^{\prime }\cup C$ and
$G_{1}^{\prime \prime }=G-G_{2}^{\prime \prime }$. Hence suppose $C$ has
length $\equiv 2\pmod 3$. Since $G_{1}$ has no cycle of $\mathcal{F}$, we
may assume that $C$ contains one of the vertices $u_{j+1}^{2},\ldots
,u_{m_{2}}^{2}$. {Let }$\ell \in \{j+1,\ldots ,m_{2}\}$ be the smallest
index {such that }$u_{\ell }^{2}\in V(C)$. {Let }$M^{\prime }=(M\cup
C)+u_{j}^{2}u_{j+1}^{2}\ldots u_{\ell }^{2}$. Clearly $M^{\prime }$ is
strong {because it belongs to }$\mathcal{B}_{r,s},$ {with }$r+s\geq 3$ and $%
s\leq 2$. {By considering the subgraphs }$G_{2}^{\prime \prime
}=(G_{2}-(M_{1}\cup M_{2}))\cup M^{\prime }$ and $G_{1}^{\prime \prime
}=G-G_{2}^{\prime \prime },$ {the pair }$(G_{1}^{\prime \prime
},G_{2}^{\prime \prime })$ {leads to a contradiction on the choice of }$%
(G_{1},G_{2}).$

Now let $G_{1}^{\prime }$ has at least two disjoint cycles $C$ and $%
C^{\prime }$ of $\mathcal{F}$. Using an argument similar to that described
in the proof of Lemma \ref{mainn}, we can obtain a pair $(G_{1}^{\prime
\prime },G_{2}^{\prime \prime })$ such that $G_{1}^{\prime \prime }$ has no
cycle of $\mathcal{F}$ and each component of $G_{2}^{\prime \prime }$
belongs to $\mathcal{F}_{0}\cup \mathcal{F}_{0,2}\cup \mathcal{F}_{2,2}\cup
\mathcal{F}_{3}\cup \mathcal{E}$, {where either }$G_{2}^{\prime \prime }$ {%
has more strong components than }$G_{2}$ or the number of cycles of $%
G_{2}^{\prime \prime }$ {belonging to }$\mathcal{F}$ is greater than the\
number of cycles of $G_{2}$ {\ belonging to}$\mathcal{F}$ or $%
s_{(G_{1}^{\prime \prime },G_{2}^{\prime \prime })}>s_{(G_{1},G_{2})}.$ {In
either case, we obtain a }contradiction. $\hfill\Box $

\bigskip

Recall that a component of $\mathcal{B}_{r,s}$ is not strong when $s\geq 3.$

\smallskip \noindent \textbf{Claim 3.} Let $M_{1}$ and $M_{2}$ be two
non-strong components of $G_{2}$ such that $M_{1}\in \mathcal{F}_{2,2}$ and $%
M_{2}\in \mathcal{B}_{r,s}$. Then there is no path $v_{0}v_{1}\ldots
v_{t+1}\;(t\geq 1)$ in $G$ such that $v_{1},\ldots ,v_{t}\in V(G_{1})$, $%
v_{0}\in M_{1}$, $v_{t+1}\in M_{2}$ and $v_{t+1}$ is\ not {the }special
vertex of $M_{2}$.\newline
\noindent \textbf{Proof of Claim 3.} Suppose, to the contrary, {\ that }%
there is a path $v_{0}v_{1}\ldots v_{t+1}\;(t\geq 1)$ in $G$ such that $%
v_{1},\ldots ,v_{t}\in V(G_{1})$, $v_{0}\in M_{1}$, $v_{t+1}\in M_{2}$ and $%
v_{t+1}$ is not\textbf{\ }special vertex of $M_{2}$. Suppose $M_{1}$ is
obtained from two cycles $C^{1}=u_{1}^{1}\ldots u_{m_{1}}^{1}u_{1}^{1}$ and $%
C^{2}=u_{1}^{2}\ldots u_{m_{2}}^{2}u_{1}^{2}$ by adding the edge $%
u_{1}^{1}u_{1}^{2}$ and let $M_{2}$ obtained from $r\geq 0$ tailed-cycle $%
C_{n_{1},\ell _{1}}$, $\ldots ,C_{n_{r},\ell _{r}}$ and $s\geq 3$ cycles $%
C_{m_{1}},\ldots ,C_{m_{s}}$, where $n_{i}\equiv 2\pmod 3$ and $m_{j}\equiv 2%
\pmod
3$ for each $i,j$, by adding a new vertex $z$ (special vertex) and attaching
$z$ to the leaf of each tailed cycle $C_{n_{i},\ell _{i}}$ and to one vertex
of each cycle $C_{m_{j}}$. Without loss of generality, that we may assume
that $v_{1}$ is adjacent to the vertex $u_{j}^{2}\in V(C_{2})$ where $j\geq
m_{2}/2$ ({by relabeling the vertices if necessary)}.

First let $v_{t+1}$ belongs to a cycle $C_{m_{i}}=w_{1}^{i}w_{2}^{i}\ldots
w_{m_{i}}^{i}w_{1}^{i}$ for some $i.$ {Without loss of generality, let }$i=1$
and {$v_{t+1}=w_{q}^{1}$.} Let $M_{2}^{\prime }$ be obtained
from $M_{2}$ by deleting the vertices of $C_{m_{1}}$, and $M_{1}^{\prime
}=(C^{1}\cup C_{m_{1}})\cup u_{1}^{1}u_{1}^{2}u_{2}^{2}\ldots
u_{j}^{2}v_{1}\ldots v_{t}w_{q}^{1}$. {In this case, consider the subgraphs }%
$G_{2}^{1}=(G_{2}-(M_{1}\cup M_{2}))\cup (M_{1}^{\prime }\cup M_{2}^{\prime
})$ and $G_{1}^{1}=G-G_{2}^{1}$ {{which we will be discussing later}. }

Now {assume that }$v_{t+1}$ {belongs to a tailed cycle }$C_{n_{i},\ell _{i}}$
for some $i$, say $i=1$. Let $C=w_{1}^{1}w_{2}^{1}\ldots
w_{n_{1}}^{1}w_{1}^{1}$ be the cycle of $C_{n_{1},\ell _{1}}$ and $%
P=y_{1}^{1}\ldots y_{\ell _{1}}^{1}$ be the tail of $C_{n_{1},\ell _{1}}$
such that $w_{1}^{1}y_{1}^{1}\in E(G)$. {Consider the two situations
depending on whether }$v_{t+1}$ {is on the cycle or the tail. If }$%
v_{t+1}\in V(C)$, say $v_{t+1}=w_{q}^{1}$, {then let }$M_{2}^{\prime }$ be
obtained from $M_{2}$ by deleting the vertices of $C_{n_{1},\ell _{1}}$, and
$M_{1}^{\prime }=(C^{1}\cup C)\cup u_{1}^{1}u_{1}^{2}u_{2}^{2}\ldots
u_{j}^{2}v_{1}\ldots v_{t}w_{q}^{i}$. {In this case, consider the subgraphs }%
$G_{2}^{1}=(G_{2}-(M_{1}\cup M_{2}))\cup (M_{1}^{\prime }\cup M_{2}^{\prime
})$ and $G_{1}^{1}=G-G_{2}^{1}.$ {If }$v_{t+1}\in V(P),$ say $%
v_{t+1}=y_{q}^{1}$, {then let }$M_{2}^{\prime }$ be obtained from $M_{2}$ by
deleting the vertices of $C_{n_{1},\ell _{1}}$, and $M_{1}^{\prime
}=(C^{1}\cup C)\cup u_{1}^{1}u_{1}^{2}u_{2}^{2}\ldots u_{j}^{2}v_{1}\ldots
v_{t}y_{q}^{1}y_{q-1}^{1}\ldots y_{1}^{1}w_{1}^{1}$. {In this case, consider
the subgraphs }$G_{2}^{1}=(G_{2}-(M_{1}\cup M_{2}))\cup (M_{1}^{\prime }\cup
M_{2}^{\prime })$ and $G_{1}^{1}=G-G_{2}^{1}$.

{Observe that in any situation, either }the number of cycles of $G_{2}^{1}$ {%
belonging to }$\mathcal{F}$ is greater than {the one of }$G_{2}$ {that are
in }$\mathcal{F}$ or $s_{(G_{1}^{1},G_{2}^{1})}>s_{(G_{1},G_{2})}$. {Now, if
}$G_{1}^{1}$ has no cycle of $\mathcal{F}$, then the pair $%
(G_{1}^{1},G_{2}^{1})$ leads to a contradiction. Otherwise, by repeating
above process we can obtain a pair $(G_{1}^{\prime },G_{2}^{\prime })$ such
that $G_{1}^{\prime}$ has no cycle of $\mathcal{F}$ and each component of $%
G_{2}^{\prime}$ belongs to $\mathcal{F}_{0}\cup \mathcal{F}_{0,2}\cup
\mathcal{F}_{2,2}\cup \mathcal{F}_{3}\cup \mathcal{E}$, {where either the
number of }strong components of $G_{2}^{\prime }$ is greater than the {one
of }$G_{2}$ or the number of cycles of $G_{2}^{\prime }$ {that are in }$%
\mathcal{F}$ is greater than the number of cycles of $G_{2}$ {belonging to }$%
\mathcal{F}$ or $s_{(G_{1}^{\prime },G_{2}^{\prime })}>s_{(G_{1},G_{2})}.$ {%
In either case, we have a contradiction and the desired claim follows. }$%
\hfill\Box $


\bigskip

\smallskip \noindent \textbf{Claim 4.} {If }$M\in \mathcal{B}_{r,s}\;${is a
non-strong component of }$G_{2}$ {\ with a special vertex }$z,$ then there
is no path $v_{0}v_{1}\ldots v_{t}v_{t+1}$ $(t\geq 1)$ in $G$ such that $%
v_{1},\ldots ,v_{t}\in V(G_{1})$, $v_{0},v_{t+1}\in V(M)-\{z\}$ and $%
v_{0},v_{t+1}$ belong to different {near cycles of }$z.$ \newline
\noindent \textbf{Proof of Claim 4.} Let $M$ be obtained from $r\geq 0$
tailed-cycle $C_{n_{1},\ell _{1}}$, $\ldots ,C_{n_{r},\ell _{r}}$ and $s\geq
3$ cycles $C_{m_{1}},\ldots ,C_{m_{s}}$, where $n_{i}\equiv 2\pmod 3$ and $%
m_{j}\equiv 2\pmod
3$ for each $i,j$, by adding a new vertex $z$ (special vertex) and attaching
$z$ to the leaf of each tailed cycle $C_{n_{i},\ell _{i}}$ and to one vertex
of each cycle $C_{m_{j}}$. {Moreover, let }$C_{m_{i}}=z_{1}^{i}z_{2}^{i}%
\ldots z_{m_{i}}^{i}z_{1}^{i}$ for {\ each }$i\in \{1,...,s\}$ and let $%
V(C_{n_{i},\ell _{i}})=\{x_{1}^{i},\ldots ,x_{n_{i}}^{i},y_{1}^{i},\ldots
,y_{\ell _{i}}^{i}\},$ where $x_{1}^{i},\ldots ,x_{n_{i}}^{i}$ {induce in
order the cycle }of $C_{n_{i},\ell _{i}}$ and $y_{1}^{i},\ldots ,y_{\ell
_{i}}^{i}$ {induce in order the }tail of $C_{n_{i},\ell _{i}}$.

Suppose, to the contrary, that there is a path $P=v_{0}v_{1}\ldots
v_{t}v_{t+1}$ in $G$ such that $v_{1},\ldots ,v_{t}\in V(G_{1}),$ $%
v_{0},v_{t+1}\in V(M)-\{z\}$ and $v_{0},v_{t+1}$ belong to different {near
cycles of }$z.$

First let $r+s=3$. Then $r=0$ and $s=3$. Assume, without loss of generality,
that $v_{0}=u_{k}^{1}$ and $v_{t+1}=u_{j}^{2}$ where $j\leq m_{1}/2$ and $%
k\leq m_{2}/2$. Let $M^{\prime }$ be obtained from $C_{1},C_{3}$ by adding
the path $Pu_{j-1}^{2}\ldots u_{1}^{2}zu_{1}^{3}$. Note that if $%
v_{t+1}=u_{1}^{2}$, then the added path will be {simply }$Pzu_{1}^{3}$. {%
Consider the subgraph }$G_{2}^{\prime }=(G_{2}-M)\cup M^{\prime }$. If $%
G-G_{2}^{\prime }$ has no cycle of $\mathcal{F}$, then the pair $%
(G-G_{2}^{\prime },G_{2}^{\prime })$ {provides a number of strong components
in }$G_{2}^{\prime }$ {greater than the one of }$G_{2},$ {contradicting our
choice of the pair }$(G_{1},G_{2}).$ {Assume now that }$G-G_{2}^{\prime }$
has exactly\textbf{\ }one cycle $C$ of $\mathcal{F}$. Then $V(C)$ meets at
least a vertex of $\{u_{j+1}^{2},\ldots ,u_{m_{2}}^{2}\} $ and let $p$ be
the largest integer that $u_{p}^{2}\in V(C) $. Let $M^{\prime \prime }$ be
obtained from $M^{\prime }\cup C$ by adding the path $u_{1}^{2}\ldots
u_{p}^{2}$. {Consider the subgraph }$G_{2}^{\prime \prime }=(G_{2}-M^{\prime
})\cup M^{\prime \prime } $. Then, {as above, }the pair $(G-G_{2}^{\prime
\prime },G_{2}^{\prime \prime })$ leads to a contradiction. {Hence we can
assume that }$G-G_{2}^{\prime }$ has at least two disjoint cycles of $%
\mathcal{F}$. Clearly each of these cycles meets at least a vertex of $%
\{u_{j+1}^{2},\ldots ,u_{m_{2}}^{2}\}$. Consider the subgraph $G^{\prime }$
of $G-G_{2}^{\prime }$ induced by the vertices of these cycles and the
vertices of $\{u_{j+1}^{2},\ldots ,u_{m_{2}}^{2}\}$, {and let }$K_{1}$ and $%
K_{2}$ be two disjoint subgraphs of {$G^{\prime }$} satisfying
the conditions of Lemma \ref{mainn}. Then the pair $(G_{2}^{\prime }\cup
K_{2},G-(G_{2}^{\prime }\cup K_{2}))$ leads to a contradiction because the
number of strong components of $G_{2}^{\prime }\cup K_{2}$ is greater than
the number of strong components of $G_{2}$.

Now let $r+s\geq 4$, {and assume that }$P$ connects two cycles $C$ and $%
C^{\prime }$ of $M$ {that are at distance one from }$z.$ Let $M^{\prime }$
be obtained from $M$ by deleting the vertices of $C\cup C^{\prime }$, and
let $M^{\prime \prime }=C+C^{\prime }+P$. {Now, if we consider the subgraphs
}$G_{2}^{\prime \prime }=(G_{2}-M)\cup (M^{\prime }\cup M^{\prime \prime })$
and $G_{1}^{\prime \prime }=G-G_{2}^{\prime \prime }$, {then one can see, as
above, that the pair }$(G_{1}^{\prime \prime },G_{2}^{\prime \prime })$
leads to a contradiction. $\hfill\Box $

\bigskip

\smallskip \noindent \textbf{Claim 5.} Let $M_{1},M_{2}\in \mathcal{E}$ be
two non-strong components of $G_{2}$ and let $z_{i}$ be the special vertex
of $M_{i}$. Then there is no path $P=v_{0}v_{1}\ldots v_{t}v_{t+1}\;(t\geq
1) $ in $G$ such that $v_{1},\ldots v_{t}\in V(G_{1})$, $v_{0}$ belongs to a
{near cycle }$C_{1}$ {of }$z_{1}$ and $v_{t+1}$ belongs to {near cycle }$%
C_{2}$ {of }$z_{2}$.\newline
\noindent \textbf{Proof of Claim 5.} {Suppose to the contrary that such a
path }$P$ exists. Let $M_{i}^{\prime }$ be obtained from $M_{i}$ by deleting
the vertices of $V(C_{i})$ for {each }$i\in \{1,2\} $ and {let }$%
M=(C_{1}\cup C_{2})+P$. {\ Consider the subgraphs }$G_{2}^{\prime
}=(G_{2}-(M_{1}\cup M_{2}))\cup (M_{1}^{\prime }\cup M_{2}^{\prime }\cup M)$
and $G_{1}^{\prime }=G-G_{2}^{\prime }$. {Since }$G_{2}^{\prime }${{\ {has
more }}}strong components than $G_{2},$ {the pair }$(G_{1}^{\prime
},G_{2}^{\prime })$ contradicts the choice of the {pair }$(G_{1},G_{2}).$ $%
\hfill\Box $

\bigskip

{Now, let }$K_{0}$ be the subgraph of $G_{2}$ that consists of all
non-strong components of $G_{2}$ and {let }$H_{0}$ be the subgraph of $G_{2}$
that consists of all strong components of $G_{2}$. By Lemma \ref{strong},
each component $M$ of $H_{0}$ has a $3$-tuple $\overrightarrow{f_{M}}$ of
RDFs of $M$ such that $\omega (\overrightarrow{f_{M}})\leq \frac{(4k+8)3n(M)%
}{6k+11}$ and all vertices of $M$ are $\overrightarrow{f_{M}}$-strong. {%
Therefore, by combining }these $3$-tuples we obtain a 3-tuple $%
\overrightarrow{f_{0}}$ of RDFs of $H_{0}$ such that $\omega (%
\overrightarrow{f_{0}})\leq \frac{(4k+8)3n(H_{0})}{6k+11}$ and all vertices
of $H_{0}$ are $\overrightarrow{f_{0}}$-strong. {Moreover, }by Lemmas \ref%
{non-strong1} and \ref{non-strong2}, each component $M$ of $K_{0}$ has a
3-tuple $\overrightarrow{g_{M}}$ of RDFs of $M$ such that $\omega (%
\overrightarrow{g_{M}})\leq \frac{(4k+8)3n(M)}{6k+11},$ and {if further }$%
M\in \mathcal{B}_{r,s}\;(s\geq 3)$, then its special vertex {as well as all
all vertices }on tailed cycle are $\overrightarrow{g_{M}}$-strong. {%
Therefore, by combining }these $3$-tuples we obtain a $3$-tuple $%
\overrightarrow{g_{0}}$ of RDFs of $K_{0}$ such that $\omega (%
\overrightarrow{g_{0}})\leq \frac{(4k+8)3n(K_{0})}{6k+11}$ and all vertices
of $K_{0}$ are $\overrightarrow{g_{0}}$-strong {except vertices on near
cycles of some }special vertex or vertices on the component in $\mathcal{F}%
_{2,2}$.

%
%

If {there is a path }$P_{1}=v_{0}v_{1}\ldots v_{t}v_{t+1}$ in $G$ such that $%
v_{1},\ldots ,v_{t}\in V(G_{1})$, $N(v_{1})\cup N(v_{t})\subseteq
V(H_{0}\cup K_{0})\cup V(P_{1})$ and both $v_{0},v_{t+1}$ belong to a
component $M\in \mathcal{F}_{2,2}$ of $K_{0}$, then we deduce from Claim 1
that both of $v_{0},v_{t+1}$ belong to {same }cycle of $M$. {Let }$M^{\prime
}=M+P$. It follows from Lemma \ref{MainLem}-(5) {and Lemmas }\ref{ear1} and %
\ref{induced} that $M^{\prime }$ has a 3-tuple $\overrightarrow{f^{\prime }}$
of RDFs of $M^{\prime }$ such that $\omega (\overrightarrow{f^{\prime }}%
)\leq 2n(M^{\prime })+2\leq \frac{(4k+8)3n(M^{\prime })}{6k+11}$ and all
vertices {of }$M^{\prime }$ are $\overrightarrow{f^{\prime }}$-strong {%
except }$v_{1}$ and $v_{t}$. {In this case, let }$K_{0}^{1}=K_{0}-M$, $%
\overrightarrow{g_{0}^{1}}=\overrightarrow{g_{0}}|_{K_{0}^{1}}$ {(the
restriction of }$\overrightarrow{g_{0}}$ on $K_{0}^{1})$ and $%
H_{0}^{1}=H_{0}\cup (M+P)=H_{0}\cup M^{\prime }.$ {Let }$\overrightarrow{%
f_{0}^{1}}$ be a 3-tuple {\ of RDFs }obtained by combining the 3-tuples $%
f_{0}$ and $\overrightarrow{f^{\prime }}$. Clearly all vertices of $%
H_{0}^{1} $ which have {\ a }neighbor outside $H_{0}^{1}\cup K_{0}^{1}$ are $%
\overrightarrow{f_{0}^{1}}$-strong. By repeating this process we obtain two {%
\ sequences }of subgraphs $K_{0}\supseteq K_{0}^{1}\supseteq \ldots
\supseteq K_{0}^{d}$ and $H_{0}\subseteq H_{0}^{1}\subseteq \ldots \subseteq
H_{0}^{d}$ so that: \textrm{(i)} {there is no path }$P=w_{0}w_{1}\ldots
w_{r}w_{r+1}$ in $G$ with $w_{1},\ldots ,w_{r}\in V(G_{1})-(\cup
_{i=1}^{d}V(P_{i}))$, $N(w_{1})\cup N(w_{r})\subseteq V(H_{0}^{d}\cup
K_{0}^{d})\cup (\cup _{i=1}^{d}V(P_{i}))$ and both $v_{0},v_{t+1}$ belong to
a component $M\in \mathcal{F}_{2,2}$ of $K_{0}^{d}$, and \textrm{(ii)} $%
H_{0}^{d}$ has a 3-tuple $\overrightarrow{f_{0}^{d}}$ such that all {its }%
vertices which have{a} neighbor outside $V(H_{0}^{d}\cup K_{0}^{d})$ are $%
\overrightarrow{f_{0}^{d}}$-strong. Let $H_{1}=H_{0}^{d}$, $K_{1}=K_{0}^{d}$%
, $\overrightarrow{g_{1}}=\overrightarrow{g_{0}}|_{K_{0}^{d}}$ and $%
\overrightarrow{f_{1}}=f_{0}^{d}$. {Observe that $\omega (%
\overrightarrow{f_{1}})\leq \frac{(4k+8)3n(H_{1})}{6k+11}$ and $\omega (%
\overrightarrow{g_{1}})\leq \frac{(4k+8)3n(K_{1})}{6k+11}$. If $%
V(G)=V(H_{1}\cup K_{1})$, then by combining 3-tuple $\overrightarrow{f_{1}}$
of $H_{1}$ and 3-tuple $\overrightarrow{g_{1}}$ of $K_{1}$, we get a 3-tuple
$\overrightarrow{h}$ of $G$ such that $\omega (\overrightarrow{h})\leq \frac{%
(4k+8)3n(G)}{6k+11}$ which will prove the theorem.} {Hence assume }that $%
V(G)\neq V(H_{1}\cup K_{1})$, {and let }$G_{2}^{1}=H_{1}\cup K_{1}$ and $%
G_{1}^{1}=G-G_{2}^{1}$.

If {there is a path }$P_{1}=v_{0}v_{1}\ldots ,v_{t}v_{t+1}$ in $G$ such that
$v_{0},v_{t+1}\in V(H_{1})$, $v_{1},\ldots ,v_{t}\in V(G_{1}^{1})$ and $%
N(v_{1})\cup N(v_{t})\subseteq V(H_{1})\cup V(P_{1})$, then let $%
H_{1}^{1}=H_{1}+P$. By Lemma \ref{ear1}, we can extend $\overrightarrow{f_{1}%
}$ to a 3-tuple $\overrightarrow{f_{1}^{1}}$ of RDFs of $H_{1}^{1}$ such
that $\omega (\overrightarrow{f_{1}^{1}})\leq \frac{(4k+8)3n(H_{1}^{1})}{%
6k+11}$, {where }all vertices of $H_{1}^{1}$ but $v_{1}$ {and }$v_{t}$ are $%
\overrightarrow{f_{1}^{1}}$-strong. {Now, if there is a path }$%
P_{2}=z_{0}z_{1}\ldots ,z_{m}z_{m+1}\;(m\geq 1)$ in $G$ such that $%
z_{0},z_{m+1}\in V(H_{1}^{1})$, $z_{1},\ldots ,z_{m}\in
V(G_{1}^{1})-V(P_{1}) $, $N(z_{1})\cup N(z_{m})\subseteq V(H_{1}^{1})\cup
V(P_{2})$, then let $H_{1}^{2}=H_{1}^{1}+P_{2}$. By Lemma \ref{ear1}, we can
extend $\overrightarrow{f_{1}^{1}}$ to a 3-tuple $\overrightarrow{f_{1}^{2}}$
of RDFs of $H_{1}^{2}$ such that $\omega (\overrightarrow{f_{1}^{2}})\leq
\frac{(4k+8)3n(H_{1}^{2})}{6k+11}$ and all new vertices but $z_{1},z_{m}$
are $\overrightarrow{f_{1}^{2}}$-strong. By repeating this process we obtain
a {sequence }of subgraphs $H_{1}\subseteq H_{1}^{1}\subseteq \ldots
\subseteq H_{1}^{q}$ so that {there is no path }$P=w_{0}w_{1}\ldots
,w_{r}w_{r+1}$ in $G$ such that $w_{0},w_{r+1}\in V(H_{1}^{q})$, $%
w_{1},\ldots ,w_{r}\in V(G_{1}^{1})-(\cup _{i=1}^{q}V(P_{i}))$ and $%
N(w_{1})\cup N(w_{r})\subseteq V(H_{1}^{q})\cup V(P)$. Moreover, $H_{1}^{q}$
has a 3-tuple $\overrightarrow{f_{1}^{q}}$ of RDFs of $H_{1}^{q}$ such that $%
\omega (\overrightarrow{f_{1}^{q}})\leq \frac{(4k+8)3n(H_{1}^{q})}{6k+11}$
and all vertices are $\overrightarrow{f_{1}^{q}}$-strong unless the vertices
which have no neighbors outside of $H_{1}^{q}\cup K_{1}$. If $%
V(G)=V(H_{1}^{q}\cup K_{1})$, {then as above, by combining
3-tuple $\overrightarrow{f_{1}^{q}}$ and 3-tuple $\overrightarrow{g_{2}}=%
\overrightarrow{g_{1}}$, the result follows.} {Hence }assume that $V(G)\neq
V(H_{1}^{q}\cup K_{1})$, {{{and let }}}$H_{2}=H_{1}^{q}$, $%
\overrightarrow{f_{2}}=\overrightarrow{f_{1}^{q}}$, $K_{2}=K_{1}$, $%
\overrightarrow{g_{2}}=\overrightarrow{g_{1}}$, $G_{2}^{2}=H_{2}\cup K_{2}$
and $G_{1}^{2}=G-G_{2}^{2}$. {In the following we will use Lemma %
\ref{12} by applying its three items, one by one (in any order), starting
with the subgraph }$G_{1}^{2}$ {and obtaining each time (when the
item occurs) a sequence of subgraphs. The last subgraph of the sequence will
be used for the next item. }

\smallskip \noindent \textbf{Case 1.} {$G_{1}^{2}$} contains a
tailed $m$-cycle $C_{m,\ell }\;({m\equiv 1\pmod 3})$, with
vertex set $\{x_{1},\ldots ,x_{m},y_{1},\ldots ,y_{\ell }\}$, such that $%
y_{\ell }$ is adjacent to some vertex $x$ of $G_{2}^{2}$ and $%
N_{G}(x_{m})\subseteq {V(G_{2}^{2})}\cup V(C_{m,\ell })$.\newline
First {assume that }$x\in V(H_{2})$. Then $x$ is $\overrightarrow{f_{2}}$%
-strong. {Let }$H_{2}^{1}$ be obtained from $H_{2}$ by adding the tailed
cycle $C_{m,\ell }$ and the edge $xy_{\ell }$. By Lemma \ref%
{tailedcycle-3p+1}, $\overrightarrow{f_{2}}$ can be extended to a 3-tuple $%
\overrightarrow{f_{2}^{1}}$ or RDFs of $H_{2}^{1}$ such that $\omega (%
\overrightarrow{f_{2}^{1}})\leq \frac{(4k+8)3n(H_{2}^{1})}{6k+11}$ and all
new vertices but $x_{m}$ are $\overrightarrow{f_{2}^{1}}$-strong. Set also $%
K_{2}^{1}=K_{2}$ and $\overrightarrow{g_{2}^{1}}=\overrightarrow{g_{1}}$.

\noindent Now {assume that }$x$ belongs to a component $M$ of $K_{2}$ such
that $M\in \mathcal{F}_{2,2}$. {Let }$M^{\prime }$ be obtained from $M$ and $%
C_{m,\ell }$ by adding the edge $xy_{\ell }$. By Lemma \ref{MainLem}{-(5)
and Lemma }\ref{tailedcycle-3p+1}, one can see that $M^{\prime }$ has a
3-tuple $f_{M^{\prime }}$ of RDFs such that $\omega (\overrightarrow{%
f_{M^{\prime }}})\leq 2n(M^{\prime })+2\leq \frac{(4k+8)3n(M^{\prime })}{%
6k+11}$ and all of its vertices but $x_{m}$ are $\overrightarrow{%
f_{M^{\prime }}}$-strong. {Let }$H_{2}^{1}=H_{2}\cup M^{\prime }$, $%
K_{2}^{1}=K_{2}-M$, $g_{1}^{2}$ {be }the restriction of $g$ on $K_{2}^{1}$
and 3-tuple $\overrightarrow{f_{2}^{1}}$ is obtained from combining $%
\overrightarrow{f_{M^{\prime }}}$ and $\overrightarrow{f_{2}}$. Note that
all vertices of $H_{2}^{1}$ which have neighbor in {$%
G_{1}^{2}-V(C_{m,\ell})$} are $\overrightarrow{f_{2}^{1}}$-strong.

\noindent Next assume that $x$ belongs to a component $M$ of $K_{2}$ such
that $M\in \mathcal{B}_{r,s}\;({s\geq 3})$ and $x$ is $%
\overrightarrow{g}$-strong. {{{Let }}}$M^{\prime }$ {{{%
be }}}obtained {{{from }}}$M$ and $C_{m,\ell }$ by
adding the edge $xy_{\ell }$ and let $H_{2}^{1}=H_{2}$ and $%
K_{2}^{1}=(K_{2}-M)\cup M^{\prime }$. By Lemma \ref{tailedcycle-3p+1}, $%
\overrightarrow{g_{2}}$ can be extended to a 3-tuple $g_{2}^{1}$ of RDFs of $%
K_{2}^{1}$ such that $\omega (\overrightarrow{g_{2}^{1}})\leq \frac{%
(4k+8)3n(K_{2}^{1})}{6k+11}$ and all newly added vertices but $x_{m}$ are $%
\overrightarrow{g_{2}^{1}}$-strong.

Finally, {assume that }$x$ belongs to a component $M$ of $K_{2}$ such that $%
M\in \mathcal{B}_{r,s}\;({s\geq 3})$ and $x$ is not $%
\overrightarrow{g}$-strong. Then $x$ belongs to a {near }cycle $C$ {from the
}special vertex of $M$. Let $M^{\prime }$ be obtained from $M$ by deleting
the vertices of $C$ and {let }$M^{\prime \prime }$ be obtained from $C$ and $%
C_{m,\ell }$ by adding the edge $xy_{\ell }$. {In this case, let }$%
H_{2}^{1}=H_{2}\cup M^{\prime \prime }$, $K_{2}^{1}=K_{2}-V(C)$ and $%
g_{1}^{2}$ is the restriction of $\overrightarrow{g_{2}}$ on $K_{2}^{1}$. By
Lemma \ref{MainLem} (items 2,3,4), $f_{2}$ can be extended to a 3-tuple of
RDFs $f_{2}^{1}$ of $H_{2}^{1}$ such that $\omega (\overrightarrow{f_{2}^{1}}%
)\leq \frac{(4k+8)3n(H_{2}^{1})}{6k+11}$ and all newly added vertices but $%
x_{m}$ are $\overrightarrow{f_{2}^{1}}$-strong.

By repeatedly applying the above argument we obtain two\ {{{
sequences }}}of subgraphs $H_{2}\supseteq H_{2}^{1}\supseteq \ldots
\supseteq H_{2}^{r_{1}}$ and {$K_{2},K_{2}^{1},\ldots
,K_{2}^{s_{1}}$} such that there is no tailed cycle $C_{m,\ell }\;({%
m\equiv 1\pmod 3})$ in $G-(H_{2}^{r_{1}}\cup K_{2}^{s_{1}})$
whose end-vertex is adjacent to a vertex of $H_{2}^{r_{1}}\cup K_{2}^{s_{1}}$%
. Let $H_{3}=H_{2}^{r_{1}}$, $K_{3}=K_{2}^{s_{1}}$, $f_{3}$ be a 3-tuple of
RDFs of $H_{3}$ such that all vertices of $H_{3}$ which have a neighbor
outside $H_{3}\cup K_{3}$ are $\overrightarrow{f_{3}}$-strong, and $g_{3}$
be a 3-tuple of RDFs of $K_{3}$ such that all newly added vertices of $K_{3}$
which have no neighbor outside $H_{3}\cup K_{3}$ are $\overrightarrow{g_{3}}$%
-strong. {Let }$G_{2}^{3}=H_{3}\cup K_{3}$ and $G_{1}^{3}=G-G_{2}^{3}$.

\smallskip \noindent \textbf{Case 2.} $G_{1}^{3}$ contains a cycle $%
C_{m}=x_{1}x_{2}\ldots x_{m}x_{1}\;({m\equiv 1\pmod 3})$ such
that $N_{G}(x_{m})\subseteq V(G_{2}^{3})\cup V(C_{m})$ and there is an edge $%
x_{1}y$ with $y\in V(G_{3}^{2})$ Applying an argument similar to that
described in Case 1, we obtained subgraphs $H_{4}$ and $K_{4}$ such that $%
H_{3}\subseteq H_{4}$, and a 3-tuple $f_{4}$ of RDFs of $H_{4}$ so that all
vertices of $H_{4}$ {having a }neighbor outside $H_{4}\cup K_{4}$ are $%
\overrightarrow{f_{4}}$-strong, and a 3-tuple $\overrightarrow{g_{4}}$ of
RDFs of $K_{4}$ so that all newly added vertices of $K_{3}$ {that have a }%
neighbor outside $H_{4}\cup K_{4}$ are $\overrightarrow{g_{4}}$-strong.
Assume that $G_{2}^{4}=H_{4}\cup K_{4}$ and $G_{1}^{4}=G-G_{2}^{4}$.

Let $\overrightarrow{h}$ be a 3-tuple defined on $G_{2}^{4}$ {\ obtained by
combining }$\overrightarrow{f_{4}}$ and $\overrightarrow{g_{4}} $.

\smallskip \noindent \textbf{Case 3.} $G_{1}^{4}$ has a path $%
P=v_{0}v_{1},\ldots ,v_{t}v_{t+1}\;(t\geq 1)$ such that $v_{0},v_{t+1}\in
G_{2}^{4}$, $v_{1},\ldots ,v_{t}\in V(G_{1}^{4})$ and $N_{G}(v_{1})\cup
N_{G}(v_{t+1})\subseteq V(G_{2}^{4})\cup V(P)$.\newline
By Claims 1,2,3,4 {and }5, at least one of the vertices $v_{0},v_{t+1}$ is $%
\overrightarrow{h}$-strong. First {assume that each of }$v_{0}$ and $v_{t+1}$
is $\overrightarrow{h}$-strong. {Let }$G_{2}^{5}=G_{2}^{4}+P$ and $%
G_{1}^{5}=G-G_{2}^{5}$. By Lemma \ref{ear1}, we can extend $\overrightarrow{h%
}$ to a 3-tuple $\overrightarrow{h^{1}}$ of RDFs of $G_{2}^{5}$ such that $%
\omega (\overrightarrow{h^{1}})\leq \frac{(4k+8)3n(G_{2}^{5})}{6k+11}$ and
all vertices of $G_{2}^{5}$ which have neighbor in $G_{1}^{5}$ are $%
\overrightarrow{h^{1}}$-strong. {Assume now, without loss of generality,
that }$v_{0}$ is $\overrightarrow{h}$-strong and $v_{t+1}$ is not $%
\overrightarrow{h}$-strong. It follows that $v_{t+1}$ is on {a near }cycle
from a special vertex of a component of $K_{4}$ or is in a component of $%
K_{4}$ {that belongs }to $\mathcal{F}_{2,2}$. If $v_{t+1}$ is on a {near }%
cycle $C$ from a special vertex, then let $G_{2}^{5}=G_{2}^{4}+P$. By Lemma %
\ref{tailedcycle-3p+1}, $h|_{G_{2}^{4}-C}$ can be extended to a 3-tuple $%
\overrightarrow{h^{1}}$ of RDFs of $G_{2}^{5}$ such that $\omega (%
\overrightarrow{h^{1}})\leq \frac{(4k+8)3n(G_{2}^{5})}{6k+11}$ and all
vertices of $G_{2}^{5}$ which have neighbors in $G_{1}^{5}-V(P)$ are $%
\overrightarrow{h^{1}}$-strong. If $v_{t_{1}}$ is in a component $M$ of $%
K_{4}$ {belonging }to $\mathcal{F}_{2,2}$, then let $G_{2}^{5}=G_{2}^{4}+P$.
By applying Lemma \ref{tailedcycle-3p+1} twice, $h|_{G_{2}^{4}-M}$ can be
extended to a 3-tuple $\overrightarrow{h^{1}}$ of RDFs of $G_{2}^{5}$ such
that $\omega (\overrightarrow{h^{1}})\leq \frac{(4k+8)3n(G_{2}^{5})}{6k+11}$
and all vertices of $G_{2}^{5}$ which have neighbors in $G_{1}^{5}-V(P)$ are
$\overrightarrow{h^{1}}$-strong.

By repeating this process we obtain a 3-tuple $\overrightarrow{h}$ of RDFs $%
G $ such that $\omega (\overrightarrow{h})\leq \frac{(4k+8)3n(G)}{6k+11}$, {%
\ implying that }$\gamma _{R}(G)\leq \frac{(4k+8)n(G)}{6k+11}$ as desired. $%
\ \ \ \hfill \Box $

\bigskip

{{{Now, the next result settling }}}Conjecture \ref{conj} is an
immediate consequence of Theorem \ref{Theorem} and the Gallai-type result $%
\gamma _{R}(G)+\partial (G)=n$ which is valid for every graph $G$ of order $%
n.$

\begin{corollary}
\emph{Let $G$ be a graph of order $n\ge 6k + 9$, minimum degree $\delta\ge 2$%
, which does not contain any induced $\{C_5,C_8,\ldots,C_{3k+2}\} $-cycles.
Then $\partial(G)\ge \frac{(2k+3)n}{6k+11}.$}
\end{corollary}

\section{Acknowledgment}
This work was supported by the National Key R \& D Program of China (Grant No. 2019YFA0706402)
and the Natural Science Foundation of Guangdong Province under grant 2018A0303130115.

\end{document}